\documentclass[a4paper]{article}
\pdfoutput=1

\usepackage{graphicx,wrapfig}
\usepackage{lipsum}
\usepackage{amsmath,amssymb,mathrsfs}
\usepackage{hyperref}
\usepackage{url}
\usepackage{mathtools}
\usepackage{changepage}
\usepackage{multirow}
\usepackage{wrapfig}
\usepackage{subfig}
\usepackage{color}
\usepackage{epsfig,enumerate,amsmath,amsfonts,amssymb,amsthm,mathrsfs,ifpdf}
\usepackage{indentfirst,relsize}
\usepackage{setspace,graphicx}
\usepackage{latexsym}
\usepackage[all]{xy}
\usepackage[usenames,dvipsnames]{pstricks}
\usepackage{pst-grad} 
\usepackage{pst-plot} 
\usepackage[margin = 2.50cm]{geometry}

\usepackage{algorithm}
\usepackage{authblk} 
\usepackage{algpseudocode}
\usepackage{color,soul}

\newcommand{\floor}[1]{\left\lfloor #1 \right\rfloor}

\newcommand{\remove}[1]{}

\newtheorem{theorem}{Theorem}
\newtheorem{lemma}[theorem]{Lemma}
\newtheorem{proposition}[theorem]{Proposition}

\newtheorem{claim}[theorem]{Claim}
\newcounter{Case}[theorem]
\newtheorem{case}[Case]{Case}

\newtheorem{definition}[theorem]{Definition}
\newtheorem{remark}{Remark}

\title{Tournament transitivity of graphs}

\author[]{Kamal Santra\footnote{kamal.7.2013@gmail.com, kamal\_1821ma04@iitp.ac.in}}
\affil[]{Department of Mathematics\\
	
	Indian Institute of Technology Patna\\
	
	Bihta, 801106, Bihar, India}
\date{}

\begin{document}

\maketitle
\begin{abstract}
	Let $G=(V, E)$ be a graph where $V$ and $E$ are the vertex and edge sets, respectively. For two disjoint subsets $A$ and $B$ of $V$, we say $A$ \textit{dominates} $B$ if every vertex of $B$ is adjacent to at least one vertex of $A$ in $G$. A vertex partition $\pi = \{V_1, V_2, \ldots, V_k\}$ of $G$ is called a \emph{transitive partition} of size $k$ if $V_i$ dominates $V_j$ for all $1\leq i<j\leq k$. A vertex partition $\pi = \{V_1, V_2, \ldots, V_k\}$ of $G$ is called a \emph{tournament transitive partition} of size $k$ if $V_i$ dominates $V_j$ for all $1\leq i<j\leq k$ and $V_j$ does not dominate $V_i$ for $i<j$. The maximum integer $k$ for which the above partition exists is called \emph{tournament transitivity} of $G$, and it is denoted by $TTr(G)$. The \textsc{Maximum Tournament Transitivity Problem} is to find a tournament transitive partition of a given graph with the maximum number of parts. In this article, we study this variation of transitive partition from a structure and algorithmic point of view. We show that the decision version of this problem is NP-complete for chordal graphs (connected), perfect elimination bipartite graphs (disconnected) and doubly chordal graphs (disconnected). On the positive side, we prove that this problem can be solved in polynomial time for trees. Furthermore, we characterize \textup{Type-I BCG} with equal transitivity and tournament transitivity and find some sufficient conditions under which the above two parameters are equal for a \textup{Type-II BCG}. Finally, we show that for \textup{Type-III BCG}, these two parameters are never equal.
\end{abstract}


\section{Introduction}

Partitioning a graph is one of the fundamental problems in graph theory. In the partitioning problem, the objective is to partition the vertex set (or edge set) into some parts with desired properties, such as independence, minimal edges across partite sets, etc. In literature, partitioning the vertex set into certain parts so that the partite sets follow particular domination relations among themselves has been studied \cite{chang1994domatic, furedi2008inequalities, haynes2019transitivity, haynes2020upper, hedetniemi2018transitivity, paul2023transitivity, santra2023transitivity, zaker2005grundy}. Let $G$ be a graph with $V(G)$ as its vertex set and $E(G)$ as its edge set. When the context is clear, $V$ and $E$ are used instead of $V(G)$ and $E(G)$. The \emph{neighbourhood} of a vertex $v\in V$ in a graph $G=(V, E)$ is the set of all adjacent vertices of $v$ and is denoted by $N_G(v)$. The \emph{degree} of a vertex $v$ in $G$, denoted as $\deg_G(v)$, is the number of edges incident to $v$. A vertex $v$ is said to \emph{dominate} itself and all its neighbouring vertices. A \emph{dominating set} of $G=(V, E)$ is a subset of vertices $D$ such that every vertex $x\in V\setminus D$ has a neighbour $y\in D$, that is, $x$ is dominated by some vertex $y$ of $D$. For two disjoint subsets $A$ and $B$ of $V$, we say $A$ \emph{dominates} $B$ if every vertex of $B$ is adjacent to at least one vertex of $A$.

There has been a lot of research on graph partitioning problems based on a domination relationship between the different sets. Cockayne and Hedetniemi introduced the concept of \emph{domatic partition} of a graph $G=(V, E)$ in 1977, in which the vertex set is partitioned into $k$ parts, say $\pi =\{V_1, V_2, \ldots, V_k\}$, such that each $V_i$ is a dominating set of $G$ \cite{cockayne1977towards}. The number representing the highest possible order of a domatic partition is called the \emph{domatic number} of G, denoted by $d(G)$. Another similar type of partitioning problem is the \emph{Grundy partition}. Christen and Selkow introduced a Grundy partition of a graph $G=(V, E)$ in 1979 \cite{CHRISTEN197949}. In the Grundy partitioning problem, the vertex set is partitioned into $k$ parts, say $\pi =\{V_1, V_2, \ldots, V_k\}$, such that each $V_i$ is an independent set and for all $1\leq i< j\leq k$, $V_i$ dominates $V_j$. The maximum order of such a partition is called the \emph{Grundy number} of $G$, denoted by $\Gamma(G)$. In 2018, J. T. Hedetniemi and S. T. Hedetniemi \cite{hedetniemi2018transitivity} introduced a transitive partition as a generalization of the Grundy partition. A \emph{transitive partition} of size $k$ is defined as a partition of the vertex set into $k$ parts, say $\pi =\{V_1,V_2, \ldots, V_k\}$, such that for all $1\leq i< j\leq k$, $V_i$ dominates $V_j$. The maximum order of such a transitive partition is called the \emph{transitivity} of $G$ and is denoted by $Tr(G)$. Recently, in 2020, Haynes et al. generalized the idea of domatic partition as well as transitive partition and introduced the concept of \emph{upper domatic partition} of a graph $G$, where the vertex set is partitioned into $k$ parts, say $\pi =\{V_1, V_2, \ldots, V_k\}$, such that for each $i, j$, with $1\leq i<j\leq k$, either $V_i$ dominates $V_j$ or $V_j$ dominates $V_i$ or both \cite{haynes2020upper}. The maximum order of such an upper domatic partition is called the \emph{upper domatic number} of $G$, denoted by $D(G)$. All these problems, domatic number \cite{chang1994domatic, zelinka1980domatically, zelinka1983k}, Grundy number \cite{effantin2017note, furedi2008inequalities, hedetniemi1982linear, zaker2005grundy, zaker2006results}, transitivity \cite{haynes2019transitivity, hedetniemi2018transitivity, paul2023transitivity, santra2023transitivity}, upper domatic number \cite{haynes2020upper, samuel2020new} have been extensively studied both from an algorithmic and structural point of view.

The concept of \emph{tournament transitive partition} was introduced by Haynes et al. in 2019 \cite{haynes2019transitivity} as an open problem. So far, no research has been done on this variation of transitivity. In this article, we study this variation of the transitivity problem from a structural and algorithmic point of view. For two disjoint subsets $A$ and $B$, we say $A$ \emph{dominates} $B$ if every vertex of $B$ is adjacent to at least one vertex of $A$. A vertex partition $\pi = \{V_1, V_2, \ldots, V_k\}$ of $G$ is called a \emph{tournament transitive partition} of size $k$ if $V_i$ dominates $V_j$ and $V_j$ does not dominate $V_i$ for all $1\leq i<j\leq k$. The maximum integer $k$ for which the above partition exists is called \emph{tournament transitivity} of $G$, and it is denoted by $TTr(G)$. A tournament transitive partition of order $TTr(G)$ is called a \emph{$TTr(G)$-partition}. The \textsc{Maximum Tournament Transitivity Problem} is to find a tournament transitive partition of a given graph with the maximum number of parts. Note that every tournament transitive partition is also a transitive partition. Therefore, for any graph $G$, $1\leq TTr(G)\leq Tr(G)\leq n$. The \textsc{Maximum Tournament Transitivity Problem} and its corresponding decision version are defined as follows.\\

\begin{center}

	\fbox{%
		\parbox{0.8\linewidth}{%
			\noindent\textsc{Maximum Tournament Transitivity Problem(MTTP)}
			
			\noindent\emph{Instance:} A graph $G=(V,E)$
			
			\noindent\emph{Solution:} A tournament transitive partition of $G$ with maximum size%
		}%
	}
	
	\vspace{0.5cm}

	\fbox{%
		\parbox{0.80\linewidth}{%
			\noindent\textsc{Maximum Tournament Transitivity Decision Problem(MTTDP)}
			
			\noindent\emph{Instance:} A graph $G=(V,E)$, integer $k$
			
			\noindent\emph{Question:} Does $G$ have a tournament transitive partition of order at least $k$?%
		}%
	}
	
\end{center}







\vspace{0.5cm}

Note that the only tournament transitive partition for the graph $K_n$ is $\{V\{K_n\}\}$. But on the other hand, according to \cite{hedetniemi2018transitivity}, we know that $Tr(K_n)=n$. This distinction motivates us to investigate this new parameter. The other similarities and differences can be found in the properties of tournament transitivity section. 

In this paper, we study this variation of transitive partition from a structure and algorithmic point of view. The main contributions are summarized below:

\begin{enumerate}
	
	\item [1.]  We show that the \textsc{Maximum Tournament Transitivity Decision Problem} is NP-complete for chordal graphs (connected), perfect elimination of bipartite graphs (disconnected), and doubly chordal graphs (disconnected).
	
	\item [2.]  We design a polynomial-time algorithm to determine the tournament transitivity of a tree.
	
	\item [3.] We characterize \textup{Type-I BCG} with equal transitivity and tournament transitivity and find some sufficient conditions under which the above two parameters are equal for a \textup{Type-II BCG}. 
	
	\item [4.] We show that for \textup{Type-III BCG}, transitivity and tournament transitivity are never equal.
	
\end{enumerate}

The rest of the paper is organized as follows: Section 2 contains basic definitions and notations that are followed throughout the article. This section also discusses the properties of tournament transitivity of graphs. Section 3 shows that the \textsc{Maximum Tournament Transitivity Decision Problem} is NP-complete in chordal graphs (Connected), perfect elimination bipartite graphs (disconnected) and doubly chordal graphs (disconnected). In Section 4, we design a polynomial-time algorithm for solving \textsc{Maximum Tournament Transitivity Problem} in trees. In section 5, we characterize \textup{Type-I BCG} with equal transitivity and tournament transitivity and find some sufficient conditions under which the above two parameters are equal for a \textup{Type-II BCG}. This section also shows that for \textup{Type-III BCG}, these two parameters are never equal. Finally, Section 6 concludes the article with some open problems.

\section{Preliminaries}

\subsection{Notation and definition}

Let $G=(V, E)$ be a graph with $V$ and $E$ as its vertex and edge sets, respectively. A graph $H=(V', E')$ is said to be a \emph{subgraph} of a graph $G=(V, E)$ if and only if $V'\subseteq V$ and $E'\subseteq E$. For a subset $S\subseteq V$, the \emph{induced subgraph} on $S$ of $G$ is defined as the subgraph of $G$ whose vertex set is $S$ and edge set consists of all of the edges in $E$ that have both endpoints in $S$, and it is denoted by $G[S]$. The \emph{complement} of a graph $G=(V,E)$ is the graph $\overline{G}=(\overline{V}, \overline{E})$, such that $\overline{V}=V$ and $\overline{E}=\{uv| uv\notin E \text{ and } u\neq v\}$. The \emph{open neighbourhood} of a vertex $x\in V$ is the set of vertices $y$ adjacent to $x$, denoted by $N_G(x)$. The \emph{closed neighborhood} of a vertex $x\in V$, denoted as $N_G[x]$, is defined by $N_G[x]=N_G(x)\cup \{x\}$. Let $S\subseteq V$, then we define $N_G(S)=\displaystyle \bigcup_{x\in S} N_G(x)$.

A subset of $S\subseteq V$ is said to be an \emph{independent set} of $G$ if every pair of vertices in $S$ are non-adjacent.  A subset of $K\subseteq V$ is said to be a \emph{clique} of $G$ if every pair of vertices in $K$ are adjacent. The cardinality of a clique of maximum size is called \emph{clique number} of $G$, and it is denoted by $\omega(G)$. A \emph{path} of length $k-1$ in a graph $G=(V, E)$, denoted by $P_k$, is a sequence of distinct vertices $v_1,\ldots,v_k$ such that $v_{i-1}v_{i} \in E$ for $2 \leq i \leq k$. Vertices $v_1$ and $v_k$ are called the \emph{end vertices} of $P_k$. For any two vertices $y, z$, we denoted a path starting with $y$ and ending at $z$ by $yPz$.


A graph is called \emph{bipartite} if its vertex set can be partitioned into two independent sets. A \emph{star} $S_t$ is the complete bipartite graph $K_{1, t}$. An edge $uv$ in a bipartite graph $G$ is called \emph{bisimplicial} if $N(u)\cup N(v)$ induces a biclique in $G$. For an edge ordering $(e_1,e_2, \ldots, e_k)$, let $S_i$ be the set of endpoints of $\{e_1,e_2, \ldots, e_i\}$ and $S_0=\emptyset$. An ordering $(e_1,e_2, \ldots, e_k)$ is a perfect edge elimination ordering for a bipartite graph $G=(V, E)$ if $G[V\setminus S_k]$ has no edges and each edge $e_i$ is a bisimplicial edge in $G[V\setminus S_{i-1}]$. A graph $G$ is a perfect elimination bipartite if and only if it admits a perfect edge elimination ordering \cite{golumbic1978perfect}.

A bipartite graph $G=(X\cup Y,E)$ is called a \textit{bipartite chain graph} if there exists an ordering of vertices of $X$ and $Y$, say $\sigma_X= (x_1,x_2, \ldots ,x_{n_1})$ and $\sigma_Y=(y_1,y_2, \ldots ,y_{n_2})$, such that $N(x_{n_1})\subseteq N(x_{n_1-1})\subseteq \ldots \subseteq N(x_2)\subseteq N(x_1)$ and $N(y_{n_2})\subseteq N(y_{n_2-1})\subseteq \ldots \subseteq N(y_2)\subseteq N(y_1)$. Such ordering of $X$ and $Y$ is called a \emph{chain ordering}, and it can be computed in linear time \cite{heggernes2007linear}. Let $G$ be a bipartite chain graph but not a complete bipartite graph. Also, let $\sigma_X= (x_1,x_2, \ldots, x_{n_1})$ and $\sigma_Y=(y_1,y_2, \ldots, y_{n_2})$ be the chain ordering of $G$ and $t$ be the maximum integer such that $G$ contains a $K_{t, t}$ as an induced subgraph. A bipartite chain graph $G$ is called $(i)$ \textup{Type-I BCG} if $x_{t+1}y_t\notin E(G)$ and $x_ty_{t+1}\notin E(G)$,  $(ii)$ \textup{Type-II BCG} if either $x_{t+1}y_t\in E(G)$ or $x_ty_{t+1}\in E(G)$ not both,  $(iii)$ \textup{Type-III BCG} if $x_{t+1}y_t\in E(G)$ and $x_ty_{t+1}\in E(G)$.

An edge between two non-consecutive vertices of a cycle is called a \emph{chord}. If every cycle in $G$ of length at least four has a chord, then $G$ is called a \emph{chordal graph}. A vertex $v\in V$ is called a \emph{simplicial vertex} of $G$ if $N_G[v]$ induces a clique in $G$. A \emph{perfect elimination ordering (PEO)} of $G$ is an ordering of the vertices, say $\sigma=(v_1, v_2,\ldots, v_n)$, such that $v_i$ is a simplicial vertex of $G_i=G[\{v_i, v_{i+1}, \ldots , v_n\}]$ for all $1\leq i\leq n$. Chordal graphs can be characterized by the existence of PEO; that is, a graph $G$ is chordal if and only if $G$ has a PEO \cite{fulkerson1965incidence}. A vertex $u\in N_G[v]$ is called a \emph{maximum neighbour} of $v$ in $G$ if $N_G[w]\subseteq N_G[u]$ for every vertex $w \in  N_G[v]$. A vertex $v$ in $G$ is called \emph{doubly simplicial} if it is simplicial and has a maximum neighbour in $G$. An ordering $\delta=(v_1, v_2,\ldots, v_n)$ of vertices of $G$ is called a \emph{doubly perfect elimination ordering (DPEO)} if $v_i$ is a doubly simplicial vertex in $G_i = G[\{v_i, v_{i+1}, \ldots , v_n\}]$ for all $1\leq i\leq n$. A graph is doubly chordal if it admits a doubly perfect elimination ordering \cite{brandstadt1998dually}.

\subsection{Properties of tournament transitivity}
In this subsection, we present some properties of tournament transitivity that motivate us to study this variation of transitivity. First, we show the following bounds for tournament transitivity.

\begin{proposition} \label{upper_bound_tournament_transitivity}
	For any graph $G$, $TTr(G)\leq min\{Tr(G), Tr(\overline{G})\}$.
\end{proposition}

\begin{proof}
	Let $\pi=\{V_1, V_2, \ldots, V_k\}$ be a $TTr$-partition of $G$. By the definition of tournament transitive partition, $\pi$ is also a transitive partition of $G$. Thus, we can say that $TTr(G)\leq Tr(G)$. Now we show that $\pi$ is also a transitive partition of $\overline{G}$. Since $\pi$ is a tournament transitive partition, for $i<j$, $V_i$ dominates $V_j$ and $V_j$ does not dominate $V_i$. Therefore, there exists a vertex $x\in V_i$ such that $xy\notin E(G)$ for all $y\in V_j$, and hence $V_i$ dominates $V_j$ in $\overline{G}$. Thus, $\pi$ is also a transitive partition of $\overline{G}$. Therefore, $TTr(G)\leq Tr(\overline{G})$. Hence, for any graph $G$, $TTr(G)\leq \min\{Tr(G), Tr(\overline{G})\}$.
\end{proof}

Next, we present two upper bounds: one in terms of the maximum and minimum degrees of a graph and the other in terms of the number of vertices.

\begin{proposition} \label{upper_bound_in_different_form_tournament_transitivity}
	For any graph $G$, $TTr(G)\leq min\{\Delta(G)+1, n-\delta(G)\}$, where $n$ is the number of vertex and $\Delta(G), \delta(G)$ are maximum and minimum degree of $G$, respectively.
\end{proposition}

\begin{proof}
	From \cite{haynes2019transitivity}, we know that $Tr(G)\leq \Delta(G)+1$, for any graph $G$. Also, from the Proposition \ref{upper_bound_tournament_transitivity}, we have $TTr(G)\leq min\{Tr(G), Tr(\overline{G})\}$. Therefore, $TTr(G)\leq min\{\Delta(G)+1, \Delta(\overline{G})+1 \}=min\{\Delta(G)+1, n-\delta(G)\}$.
\end{proof}

\begin{proposition} \label{upper_bound_in_another_form_tournament_transitivity}
	For any graph $G$, $1\leq TTr(G)\leq \floor{\frac{n+1}{2}}$, where $n$ is the number of vertex of $G$.
\end{proposition}
\begin{proof}
	Taking all vertices in one set will produce a tournament transitive partition of size $1$. So, $TTr(G)\geq 1$. Now, let $\pi=\{V_1, V_2, \ldots, V_k\}$ be a $TTr$-partition of $G$. Consider $x\in V_k$. As $\pi$ $V_i$ dominates $V_k$ in $G$ for all $1\leq i<k$, $deg_G(x)\geq k-1$. Also, we know that $V_k$ does not dominate $V_j$ for all $1\leq j<k$. Therefore, $deg_{\overline{G}}(x)\geq k-1$. So, $deg_G(x)+deg_{\overline{G}}(x)\geq 2k-1$, which implies $n-1\geq deg_G(x)+deg_{\overline{G}}(x)\geq 2k-1$. Hence, $k=TTr(G)\leq \floor{\frac{n+1}{2}}$.
\end{proof}

\begin{figure}[htbp!]
	\centering
	\includegraphics[scale=0.60]{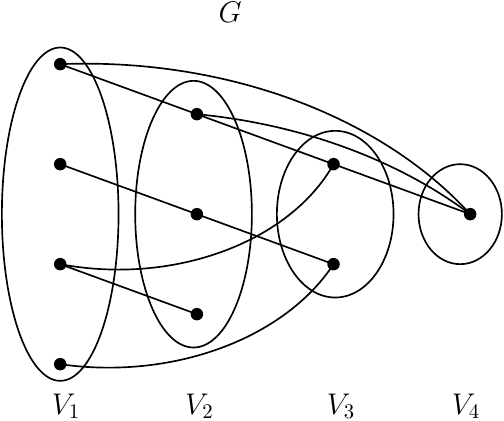}
	\caption{$\pi=\{V_1, V_2, V_3, V_4\}$ is a tournament transitive partition of $G$, but  $\pi'=\{V_1, V_2\cup V_3, V_4\}$ is not a tournament transitive partition of $G$}
	\label{fig:Arbitrarily_union}
\end{figure}

In the transitive partition $\pi=\{V_1, V_2, \ldots, V_k\}$ of a graph $G$, if we union any two sets of a transitive partition of order $k$, we create a transitive partition of order $k-1$ \cite{hedetniemi2018transitivity}. However, for the tournament transitive partition, this result is not true for arbitrary values $i$ and $j$ with $i<j$ (see the Figure\ref{fig:Arbitrarily_union}).

But if we take into account the first two sets of the partition, the result is always true. Hence, we have the following proposition.

\begin{proposition}
	If $\pi=\{V_1, V_2, \ldots, V_k\}$ is a tournament transitive partition of a graph $G$ of size $k$, then the partition $\pi'=\{V_1\cup V_2, V_3, \ldots, V_k\}$ is a tournament transitive partition of size $k-1$.
\end{proposition}

Given this, we have the following interpolation result.

\begin{proposition}\label{any_size_ttr}
	Let $G$ be a graph and $k$ be the order of a tournament transitive partition of $G$, then for every $j$, $1\leq j\leq k$, $G$ has a tournament transitive partition of order $j$.
\end{proposition}

\begin{proposition}\label{ttr_last_two_sets_size}
	Let $G$ be a connected graph and $TTr(G)=k$, $k\geq 3$. Then there exists a tournament transitive partition of $G$ of size $k$, say $\pi=\{V_1, V_2, \ldots, V_k\}$ of $G$, such that $|V_k|=1$, $|V_{k-1}|=2$. Moreover, $V_{k-1}=\{x, y\}$ and $ V_k=\{z\}$ such that $xz\in E(G)$ and $yz\notin E(G)$.
\end{proposition}

Next, we find the tournament transitivity of complete graphs, paths, cycles and complete bipartite graphs.
\begin{proposition} \label{Complete_graph_ttr}
	For the complete graph $K_n$ of $n$ vertices $TTr(K_n)=1$.
\end{proposition}
\begin{proof}
	The partition $\{V(K_n)\}$ is a tournament transitive partition of $K_n$ of size $1$. Moreover, from the Proposition \ref{upper_bound_in_different_form_tournament_transitivity}, we have $TTr(G)\leq min\{\Delta(K_n)+1, n-\delta(K_n)\}= min\{n, 1\}=1$. Therefore, $TTr(K_n)=1$.
\end{proof}

\begin{proposition}\label{Path_ttr}
	Let $P_n$ be a path of $n$ vertices, and then the tournament transitivity of $P_n$ is given as follows:
	
	$TTr(P_n) = \begin{cases}
		1 & n=1, 2 \\
		2 & n=3, 4 \\
		3 & n\geq 5
	\end{cases}$
\end{proposition}

\begin{proof}
	For $n=1, 2$, as $P_1$ and $P_2$ are $K_1$ and $K_2$, respectively, form the Proposition \ref{Complete_graph_ttr}, we have $TTr(P_n)=1$ for $n=1, 2$. Let us consider $n=3, 4$. According to Proposition \ref{upper_bound_in_another_form_tournament_transitivity}, we have $TTr(P_n)\leq 2$. Also, from the Figure\ref{fig:TTR_Path_example}, it is clear that $TTr(P_n)\geq 2$. Hence, $TTr(P_n)=2$ for $n=2, 3$. 
	
	\begin{figure}[htbp!]
		\centering
		\includegraphics[scale=0.65]{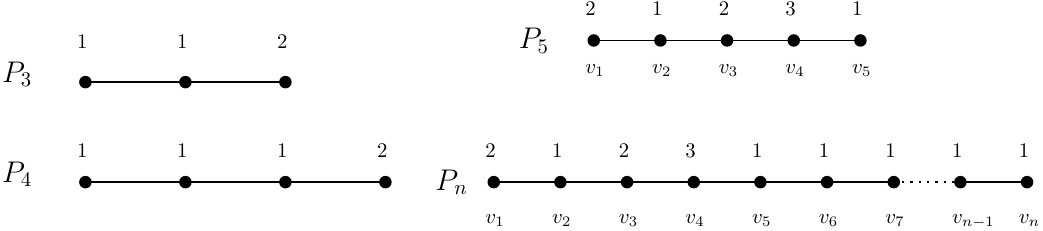}
		\caption{Tournament transitive partition of $P_n$, $n\geq 3$}
		\label{fig:TTR_Path_example}
	\end{figure} 
	
	Assume $n=5$. We show that $TTr(P_5)=3$. From the Proposition\ref{upper_bound_tournament_transitivity}, we have $TTr(P_5)\leq 3$. Let $\{v_1, v_2, v_3, v_4, v_5\}$ be the vertices of $P_5$, where $v_1, v_5$ are its end vertices. Consider a vertex partition $\pi=\{V_1,V_2, V_3\}$ as follows: $V_1=\{v_2, v_5\}$, $V_2=\{v_1, v_3\}$, and $V_3=\{v_4\}$ (see Figure \ref{fig:TTR_Path_example}). We show that $\pi$ is a tournament transitive partition. Clearly, $V_i$ dominates $V_j$, and $V_j$ does not dominate $V_i$, for all $1\leq i<j\leq j$. Therefore, $TTr(P_5)=3$. Finally, consider $n>5$. Also, let $\{v_1, v_2, v_3, v_4, v_5, V_6, \ldots, v_n\}$ be the path $P_n$. Consider a vertex partition $\pi=\{V_1,V_2, V_3\}$ as follows: $\{v_2, v_5\}\subseteq V_1$, $V_2=\{v_1, v_3\}$, and $V_3=\{v_4\}$ and all the other vertices in $V_1$ (see Figure \ref{fig:TTR_Path_example}). Similarly, we can show that $TTr(P_n)=3$ for all $n>5$. Hence, we have $TTr(P_n)=3$, $n\geq 5$.
\end{proof}

\begin{proposition}\label{Cycle_ttr}
	Let $C_n$ be a path of $n$ vertices, and then the tournament transitivity of $P_n$ is given as follows:
	
	$TTr(C_n) = \begin{cases}
		1 & n=3  \\
		2 & n=4, 5 \\
		3 & n\geq 6
	\end{cases}$
\end{proposition}

\begin{proof}
	Note that $TTr(C_n)\leq 3$, for all $n\geq 3$. We show that $TTr(C_5)=2$; for the other cases, we can easily prove the statements. If possible, assume $TTr(C_5)=3$ and $\{V_1, V_2, V_3\}$ be a tournament transitive partition of $C_5$. Also, assume $C_5=\{v_1, v_2, v_3, v_4, v_5\}$. Without loss of generality, assume $v_1\in V_3$. As $V_1, V_2$ dominate $V_3$, each of $V_1, V_2$ contains exactly one vertex from $\{v_2, v_5\}$. Again, $V_3$ does not dominate $V_1, V_2$; there exists a vertex from $\{v_3, v_4\}$ that must be in $V_2$, and a vertex from $\{v_3, v_4\}$ must be in $V_1$. In this situation, $V_1$ dominates $V_2$, and $V_2$ dominates $V_1$, a contradiction as $\pi$ is a tournament transitive partition. Hence, $TTr(C_5)=2$.
\end{proof}

\begin{proposition} \label{Complete_bipartite_graph_ttr}
	Let $G$ be a complete bipartite graph $K_{m,n}$ with $m+n$ vertices. Then tournament transitivity of $K_{m, n}$ is given as follows:
	
	$TTr(K_{m, n}) = \begin{cases}
		1 & m=n=1  \\
		2 & \text{otherwise}
	\end{cases}$
\end{proposition}
\begin{proof}
	
	Assume $m=1$ and $n=1$. Clearly, $K_{m, n}$ is a $P_2$ and from Proposition \ref{Complete_graph_ttr}, we have $TTr(K_{m, n})=1$.
	
	Next, we show that $TTr(K_{m,n})=2$ if either $m\neq 1$ or $n\neq 1$. Let the vertex set of $K_{m,n}$ be $V=X\cup Y$ such that $|X|=m, |Y|=n$ and $m\geq 2$. Consider a vertex partition $\pi=\{V_1, V_2\}$, where $V_1=(X\setminus\{x\})\cup Y$, $V_2=\{x\}$ for some $x\in X$. Since $m\geq 2$, $\pi$ is a tournament transitive partition of $K_{m,n}$ and hence $TTr(K_{m, n})\geq 2$. To prove $TTr(K_{m, n})=2$, we show that $TTr(K_{m, n})$ cannot be more than $2$. If possible, let $TTr(K_{m, n})=k\geq 3$ and $\pi=\{V_1, V_2, \ldots, V_k\}$ be a tournament transitive partition of $K_{m, n}$ of size $k$. We show that $V_i$, for $3\leq i\leq k$, cannot contain any vertex from $X$ or $Y$. If $V_i$ contains a vertex, say $x$, from $X$, there must exists $y_1\in V_1$ and $y_2\in V_2$ such that $xy_1, xy_2\in E(G)$. As $G$ is a complete bipartite graph and $V_i$ does not dominate $V_1$ and $V_2$, there must be $x_1\in V_1$ and $x_2\in V_2$. So, $V_1$ and $V_2$ contains vertices from $X$ and $Y$ and hence $V_1$ dominates $V_2$ and $V_2$ dominates $V_1$, which is a contradiction as $\pi$ is a tournament transitive partition of $K_{m,n}$. On the other hand, assume $V_i$ contains a vertex, say $y$ from $Y$. Similarly, we can show a contradiction. Therefore, $TTr(K_{m, n})=2$, when $m\geq 2$. With similar arguments, we can show that $TTr(K_{m, n})=2$, when $n\geq 2$. Hence, $TTr(K_{m, n})=2$, if either $m\neq 1$ or $n\neq 1$.
\end{proof}

For the transitivity, if $H$ is a subgraph of $G$, then $Tr(H)\leq Tr(G)$ \cite{hedetniemi2018transitivity}. But for tournament transitivity, this is not true. For example, let us consider $G=C_5$ and $H=P_5$. Clearly, $H$ is a subgraph of $G$, and $TTr(C_5)=2$ (by Proposition \ref{Cycle_ttr}), $TTr(P_5)=3$ (by Proposition \ref{Path_ttr}). Another example is $G=K_4$ and $H=C_4$, where $TTr(K_4)=1$ and $TTr(C_4)=2$.

In the above examples, neither of the subgraphs is an induced one. However, the subgraph containment relation holds for tournament transitive partition if we consider an induced subgraph instead of a subgraph.

\begin{proposition}\label{subgraph_containment_relation}
	If H is an induced subgraph of a graph G, then $TTr(H)\leq TTr(G)$.
\end{proposition}

\begin{proof}
	Let $\pi=\{V_1, V_2, \ldots, V_k\}$ be a $TTr$-transitive partition of $H$. Consider $\pi'=\{V_1', V_2', \ldots, V_k'\}$ a vertex partition of $G$ such that $V_1'= (V(G)\setminus V(H))\cup V_1$ and $V_p'=V_p$ for all $2\leq p\leq k$. As $H$ is an induced subgraph, $V_i'$ dominates $V_j'$, but $V_j'$ does not dominate $V_i'$ in $\pi'$. So, $\pi'$ is a tournament transitive partition of size $k$. Therefore, $k=TTr(H)\leq TTr(G)$.
\end{proof}

For the case of transitivity, if $G$ is a disconnected graph with connected component $C_1, C_2, \ldots, C_k$, then $Tr(G)= \max\{Tr(C_i) | 1\leq i\leq k\}$ \cite{haynes2019transitivity}. But we show that for tournament transitivity, $TTr(G)\geq \max\{TTr(C_i) | 1\leq i\leq k\}$ and the difference $TTr(G)-\max\{TTr(C_i) | 1\leq i\leq k\}$ can be arbitrarily large. First, we prove that if $G$ is a disjoint union of $t(\leq n)$ number of complete graphs of $n$ vertices, then $TTr(G)=t$.

\begin{lemma}\label{disjoint_union_ttr_example}
	Let $G$ be the disjoint union of $t(\leq n)$ number of a complete graph of $n$ vertices. Then $TTr(G)=t$.
\end{lemma}
\begin{proof}
	Let $G=K_n \cup \ldots \cup K_n$ ($t$ times). Let $\pi^j=\{V_1^j, \ldots, V_t^j, V_{t+1}^j, \ldots, V_n^j\}$ be a transitive partition of size $n$ of the $j$-th component of $G$, and also we denote the $j$-th component of $G$ as $K_n^j$. Let us consider $\pi=\{V_1, V_2, \ldots, V_t\}$ be a vertex partition of $G$ of size $t$ as follows: $$V_1=(\bigcup_{j=1}^t (V_{t+1}^j\cup V_{t+2}^j\cup \ldots \cup V_n^j))\cup (V_1^1\cup V_1^2\cup \ldots \cup V_{1}^{t-1}\cup (V_1^t\cup V_2^t\cup \ldots \cup V_t^t)),$$ $$V_j=(V_j^1\cup V_j^2\cup \ldots \cup V_{j}^{t-j}\cup (V_j^{t-j+1}\cup V_{j+1}^{t-j+1}\cup \ldots \cup V_{t}^{t-j+1}))$$ for all $2\leq j\leq t$ (See Figure \ref{fig:ttr_of_disjoint_union_complete_graphs} for reference). 
	
	\begin{figure}[htbp!]
		\centering
		\includegraphics[scale=0.65]{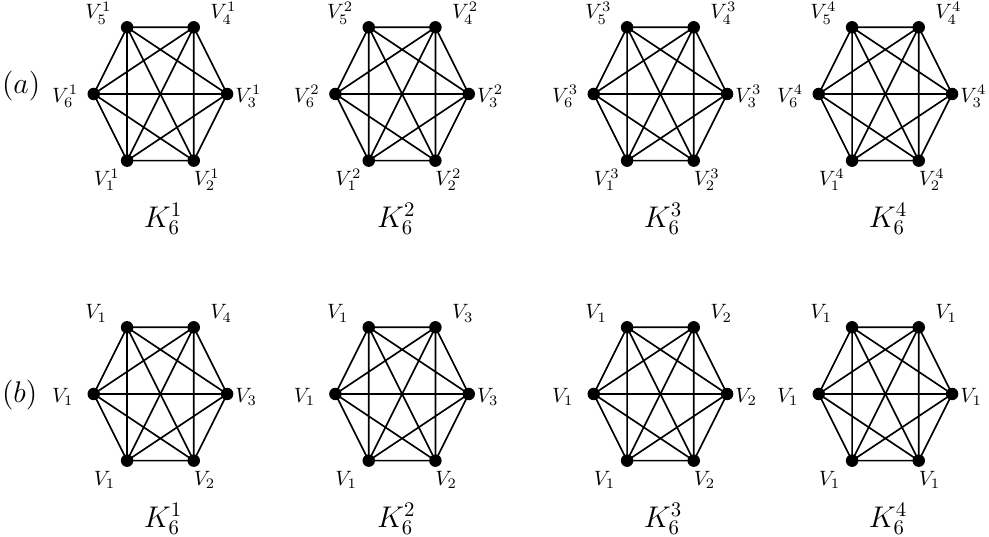}
		\caption{$(a)$ represents partitions $\pi^j$ and $(b)$ represents the partition $\pi$}
		\label{fig:ttr_of_disjoint_union_complete_graphs}
	\end{figure}

	Now we show that $\pi$ is a tournament transitive partition of size $t$. Let us consider $V_i$ and $V_j$, for $i<j$. Clearly,  $V_i\supseteq (V_i^1\cup V_i^2\cup \ldots \cup V_{i}^{t-i}\cup (V_{i}^{t-i+1}\cup V_{i+1}^{t-i+1}\cup \ldots \cup V_{t}^{t-i+1}))$ and $V_j=(V_j^1\cup V_j^2\cup \ldots \cup V_{j}^{t-j}\cup (V_j^{t-j+1}\cup V_{j+1}^{t-j+1}\cup \ldots \cup V_{t}^{t-j+1}))$. Since $\pi^p$ is a transitive partition of $K_n^p$, $V_i^p$ dominates $V_j^p$ for all $1\leq p\leq t-j$. Again, as $i<j$, which implies $t-i\geq t-j+1$ and $V_i^{t-j+1}\subseteq V_i$. So, $V_i^{t-j+1}$ dominates $(V_j^{t-j+1}\cup V_{j+1}^{t-j+1}\cup \ldots \cup V_{t}^{t-j+1})$ as $K_n^{t-j+1}$ is a complete graph. Therefore, $V_i$ dominates $V_j$ for all $1\leq i<j \leq t$. However, $V_j$ does not dominate $V_i$, as $V_i$ contains vertices from $K_n^{t-i}$, whereas $V_j$ does not contain any vertex from $K_n^{t-i}$. So, $\pi$ is a tournament transitive partition of $G$, and hence $TTr(G)\geq t$.

	To prove $TTr(G)=t$, we show that $TTr(G)$ cannot be more than $t$. If possible, let $TTr(G)=k\geq t+1$. By the Proposition \ref{any_size_ttr}, we have a tournament transitive partition of $G$ of size $t+1$, say $\pi=\{V_1, V_2, \ldots, V_{t+1}\}$. Without loss of generality, assume $V_{t+1}$ contains a vertex from $K_n^1$. As $\pi$ is a tournament transitive partition, $V_{t+1}$ does not dominate $V_{t}$, which implies there exists a vertex $x\in V_t$ and $x$ does not belong to $K_n^1$. Assume $x\in K_n^2$. Again, as $V_{t}$ does not dominate $V_{t-1}$ and $K_n^1\cap V_{t-1}, K_n^2\cap V_{t-1}$ are non-empty, there exists a vertex $y\in V_{t-1}$, and $y$ must be a vertex that does not belong to $K_n^1\cup K_n^2$. Assume $y\in K_n^3$. If we go like this, we end up with a situation where we need a vertex in $V_1$ that does not belong to $K_n^1\cup K_n^2\cup \ldots \cup K_n^t$, a contradiction. As a result, $\pi$ cannot be a tournament transitive partition of $G$, and hence $TTr(G)\leq t$. Therefore, $TTr(G)=t$.
\end{proof}

\begin{proposition}\label{disjoint_union_ttr}
	Let $G$ is a disconnected graph with connected component $C_1, C_2, \ldots, C_k$, then $TTr(G)\geq \max\{TTr(C_i) | 1\leq i\leq k\}$ and the difference $TTr(G)-\max\{TTr(C_i) | 1\leq i\leq k\}$ can be arbitrarily large.
\end{proposition}

\begin{proof}
	Let $C_i$ be the component such that $TTr(C_i)=\max\{TTr(C_i) | 1\leq i\leq k\}$. Also, let $\pi=\{V_1, V_2, \ldots, V_{k}\}$ be a $TTr$-partition of $C_i$. Now if we consider a vertex partition $\pi'=\{V_1', V_2', \ldots, V_k'\}$, where $V_1'=V_1\displaystyle \bigcup_{p=1, p\neq i}^t V(C_p)$ and $V_i'=V_i$, for all $2\leq i\leq k$. Clearly, $\pi'$ is a tournament transitive partition of $G$ of size $k$. Therefore, $TTr(G)\geq \max\{Tr(C_i) | 1\leq i\leq k\}$. 
	
	Now let us consider $G$ to be a graph, which is the disjoint union of $t (\leq n)$ number of complete graphs of $n$ vertices. Then by the Lemma \ref{disjoint_union_ttr_example} and Proposition \ref{Complete_graph_ttr}, we have $TTr(G)=t$ and $TTr(C_i)=1$, for all $i$. Therefore, the difference $TTr(G)-\max\{TTr(C_i) | 1\leq i\leq k\}=t-1$ arbitrarily large for $t$.
\end{proof}

Next, we characterize connected graphs $G$ with $TTr(G)=1$, and we give a necessary and sufficient condition for $TTr(G)\geq 2$.

\begin{proposition}\label{characterization_TTr(G)=1}
	For any connected graph $G$, $TTr(G)=1$ if and only if $G=K_n$.	
\end{proposition}

\begin{proof}
	From the Proposition\ref{Complete_graph_ttr}, it is clear that $TTr(G)=1$ if $G=K_n$. For the converse part, assume $TTr(G)=1$ and $G$ connected graph with $n\geq 3$. If $G$ is not a complete graph, there exists $x, y\in V(G)$ such that $xy\notin E(G)$. Since $G$ is connected, there exists a path connecting $x$ and $y$ with a length of at least $3$. Taking $y\in V_2$, $x\in V_1$, and other vertices are in $V_1$ form a tournament transitive partition of size $2$. As a result, $G$ must be a complete graph $K_n$.
\end{proof}

\begin{proposition}\label{characterization_TTr(G)_more_than_one}
	For any graph $G$, $TTr(G)\geq 2$ if and only if $G$ contains $P_3$ as an induced subgraph.	
\end{proposition}

\begin{proof}
	Let us assume $G$ contains $P_3$ as an induced subgraph. Then from Proposition\ref{Path_ttr} and Proposition\ref{subgraph_containment_relation}, we have $TTr(G)=2$.  Conversely, assume $TTr(G)\geq 2$. As $V_1$ dominates $V_2$ and $V_2$ does not dominate $V_1$, also $G$ is a connected graph, there must be a induce $P_3$ in $G$.
\end{proof}

\begin{remark}
	From the definition of a tournament transitive partition, it is clear that it is also a transitive partition of $G$. But the difference $Tr(G)-TTr(G)$ can be arbitrarily large. From the Proposition \ref{Complete_graph_ttr}, we have $TTr(K_n)=1$ and from \cite{haynes2019transitivity}, we know that $Tr(K_n)=n$, hence the difference $Tr(K_n)-TTr(K_n)=n-1$ is arbitrarily large for $n$. Moreover, we know $Tr(K_{m,n})=n+1$, when $n<m$. But from Proposition \ref{Complete_bipartite_graph_ttr}, we have $TTr(K_{m,n})=2$, when either $m\neq 1$ or $n\neq 1$. So, the difference $Tr(K_{m, n})-TTr(K_{m, n})=n-1$ is arbitrarily large for $n$.
\end{remark}

\section{Computational complexity}
\begin{figure}[htbp!]
	\centering
	\includegraphics[scale=0.65]{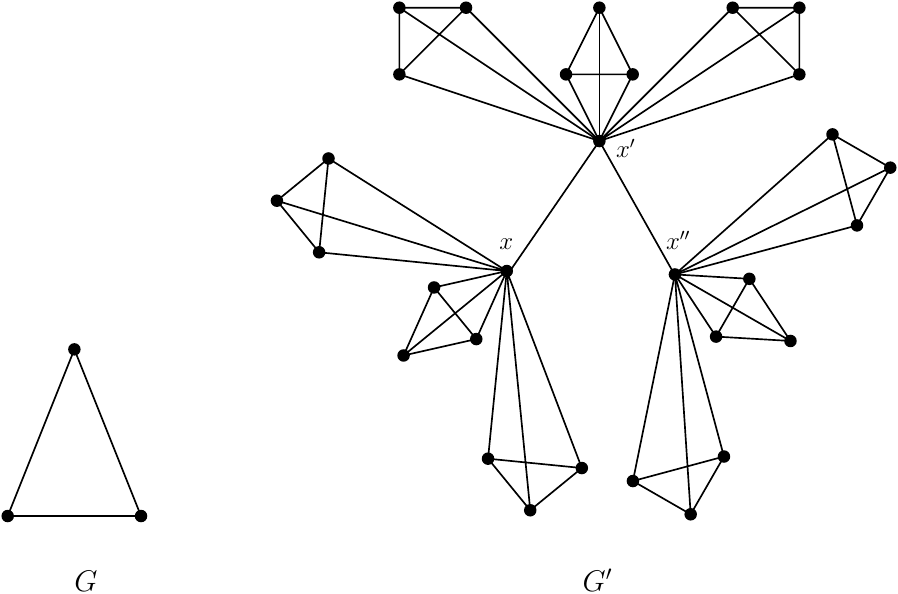}
	\caption{Construction of $G'$ from $G$}
	\label{fig:tournament_transitive_chordal_np}
\end{figure}

In this section, we present NP-completeness results for \textsc{Maximum Tournament Transitivity Decision Problem}. Given a vertex partition $\pi=\{V_1, V_2, \ldots, V_k\}$ of a graph, we can verify in polynomial time whether $\pi$ is a tournament transitive partition of that graph or not. Hence, the \textsc{Maximum Tournament Transitivity Decision Problem} is in NP. To prove that this problem is NP-hard, we show a polynomial time reduction from the \textsc{Maximum Transitivity Decision Problem} in general graphs, which is known to be NP-complete \cite{hedetniemi2004iterated}. The reduction is as follows: given an instance of the \textsc{Maximum Transitivity Decision Problem}, that is, a graph $G=(V, E)$ and an integer $k$, let us consider $3(\Delta(G)+1)$ copies of $G$. Furthermore, we take another three vertices, say $x, x'$ and $x''$. Now we connect every vertices of first $\Delta(G)+1$ copies of $G$ to $x$, second $\Delta(G)+1$ copies of $G$ to $x'$  and final $\Delta(G)+1$ copies of $G$ to $x''$. Finally, add two edges $xx'$ and $x'x''$. Clearly, the resultant graph  $G'=(V', E')$ is a chordal graph if $G$ is a chordal graph and having $(3n(\Delta(G)+1)+3)$ vertices and $(3m(\Delta(G)+1)+2)$ edges. The construction of $G'$ is illustrated in Figure \ref{fig:tournament_transitive_chordal_np}. 

Let $\mathcal{G}$ be the disjoint union of $\Delta(G)+1$ copies of $G$, and also we denote the $i$-th component of $\mathcal{G}$ as $G^i$. Next, we prove that $Tr(G)=k$ if and only if $TTr(G)=k$, as shown in the following lemma.

\begin{lemma}\label{ttr_np_lemma_1}
	Let $\mathcal{G}$ be the disjoint union of \textup{($\Delta(G)+1$)}-copies of $G$, then $Tr(G)=k$ if and only if $TTr(\mathcal{G})=k$.	
\end{lemma}
\begin{proof}
	Let $G$ be a graph with $Tr(G) = k$. Also, let $\mathcal{G}=G^1\cup G^2\cup \ldots \cup G^t$, where each $G^i$ is a graph $G$ and $t=\Delta(G)+1$. Let $\pi^j=\{V_1^j, V_2^j, \ldots, V_k^j\}$ be a transitive partition of size $k$ of the $j$-th component of $\mathcal{G}$. Let us consider $\pi=\{V_1, V_2, \ldots, V_k\}$ be a vertex partition of $\mathcal{G}$ of size $k$ as follows: $$V_1=(\displaystyle \bigcup_{j=k+1}^t V(G^j))\cup (V_1^1\cup V_1^2\cup \ldots \cup V_{1}^{k-1}\cup (V_1^k\cup V_2^k\cup \ldots \cup V_k^k)),$$ $$V_j=(V_j^1\cup V_j^2\cup \ldots \cup V_{j}^{k-j}\cup (V_j^{k-j+1}\cup V_{j+1}^{k-j+1}\cup \ldots \cup V_{k}^{k-j+1}))$$ for all $2\leq j\leq k$.

	Now we show that $\pi$ is a tournament transitive partition of size $k$. Let us consider $V_i$ and $V_j$ for $i<j$. Clearly,  $V_i\supseteq (V_i^1\cup V_i^2\cup \ldots \cup V_{i}^{k-i}\cup (V_{i}^{k-i+1}\cup V_{i+1}^{k-i+1}\cup \ldots \cup V_{k}^{k-i+1}))$ and $V_j=(V_j^1\cup V_j^2\cup \ldots \cup V_{j}^{k-j}\cup (V_j^{k-j+1}\cup V_{j+1}^{k-j+1}\cup \ldots \cup V_{k}^{k-j+1}))$. Since $\pi^p$ is a transitive partition of $G^p$, $V_i^p$ dominates $V_j^p$ for all $1\leq p\leq k-j$. Again, as $i<j$, which implies $k-i\geq k-j+1$ and $V_i^{k-j+1}\subseteq V_i$. So, $V_i^{t-j+1}$ dominates $(V_j^{k-j+1}\cup V_{j+1}^{t-j+1}\cup \ldots \cup V_{t}^{t-j+1})$ as $\pi^{k-j+1}$ is a transitive partition of $G^{k-j+1}$. Therefore, $V_i$ dominates $V_j$ for all $1\leq i<j \leq t$. However, $V_j$ does not dominate $V_i$ because $V_i$ contains vertices from $G^{k-i}$, but $V_j$ does not contain any vertex from $G^{k-i}$. So, $\pi$ is a tournament transitive partition of $\mathcal{G}$, and hence $TTr(\mathcal{G})\geq k$. To prove $TTr(\mathcal{G})=k$, we show that $TTr(\mathcal{G})$ cannot be more than $k$. By the Proposition \ref{upper_bound_tournament_transitivity}, we know that $TTr(\mathcal{G})\leq Tr(\mathcal{G})$, and also from \cite{haynes2019transitivity}, we have $Tr(\mathcal{G})=\max\{Tr(G^i), 1\leq i\leq t\}=k$. Therefore, $TTr(\mathcal{G})=k$.

	Conversely, assume $TTr(\mathcal{G})=k$. Then from \cite{haynes2019transitivity} and by the Proposition \ref{upper_bound_tournament_transitivity}, we have $k=TTr(\mathcal{G})\leq Tr(\mathcal{G})=\max\{Tr(G^i), 1\leq i\leq t\}=Tr(G)$. Clearly, $Tr(G)$ cannot be more than $k$; otherwise, $TTr(\mathcal{G})$ is more than $k$. Therefore, $Tr(G)=k$ if and only if $TTr(\mathcal{G})=k$.
\end{proof}

In the next claim, we show that $G$ has a transitive partition of size $k$ if and only if $G'$ has a tournament transitive partition of size $k+2$.

\begin{claim}
	The graph $G$ has a transitive partition of size $k$ if and only if $G'$ has a tournament transitive partition of size $k+2$.
\end{claim}

\begin{proof}
	Let $G$ be a graph with a transitive partition of size $k$. Also, let $\mathcal{G}_i$ be the disjoint union of $i$-th set of $(\Delta(G)+1)$ copies of $G$, for all $1\leq i\leq 3$. By the Lemma \ref{ttr_np_lemma_1}, we know that $TTr(\mathcal{G}_i)\geq k$, for all $i$. Now, if $\pi^i=\{V_1^i, V_2^i, \ldots, V_k^i\}$ is a transitive partition of $\mathcal{G}_i$ of size $k$, we can construct a vertex partition, say $\pi'=\{V'_1, V'_2, \ldots, V'_{k+2}\}$ of $G'$, as follows: $V'_i=V_i^1\cup V_i^2\cup V_i^3$ for all $1\leq i \leq k$ and $V'_{k+1}=\{x', x''\}, V'_{k+2}=\{x\}$. We show that $\pi'$ is a tournament transitive partition of $G'$. For any pair of sets $V'_i$ and $V'_j$ with $1\leq i< j\leq k$, $V'_i$ dominates $V'_j$ in $G'$ as $V_i^p$ dominates $V_j^p$ in $\mathcal{G}_p$ and $V'_j$ does not dominate $V'_i$ in $G'$ as $V_j^p$ does not dominate $V_i^p$ in $\mathcal{G}_p$. Also, from the construction of $G'$, it is clear that $V_i$ dominates $V_j$ but $V_j$ does not dominate $V_i$, for all $1\leq i< j\leq k+2$. Therefore, $\pi'$ is a tournament transitive partition of $G'$ of size $k+2$.
	
	Conversely, let $\pi=\{V_1,V_2,\ldots, V_{k+2}\}$ be a tournament transitive partition of $G'$ of size $k+2$. So, $Tr(G')\geq k+2$. Now consider the graph $\mathcal{G}_1'$, where $V(\mathcal{G}_1')=V(\mathcal{G}_1)\cup \{x\}$ and $E(\mathcal{G}_1')=E(\mathcal{G}_1)\cup \{xy|y\in V(\mathcal{G}_1)\}$. From \cite{hedetniemi2018transitivity}, we know that for a graph $H$ and a vertex $v$, $Tr(H)-1\leq Tr(H-v)\leq Tr(H)$, and if $H$ is a disjoint union of $C_1, C_2, \ldots, C_r$, then $Tr(H)=\max\{Tr(C_i), 1\leq i\leq r\}$. So, $Tr(G')-1\leq Tr(G'-x')=Tr(\mathcal{G}_1')$, which implies $k+1\leq Tr(G')-1\leq Tr(\mathcal{G}_1')$. Furthermore, $Tr(\mathcal{G}_1')-1\leq Tr(\mathcal{G}_1'-x)=Tr(G)$. Hence, $Tr(G)\geq k$, and $G$ has a transitive partition of size $k$.
\end{proof}

Therefore, we have the following main theorem for this section.

\begin{theorem}
	The \textsc{Maximum Transitivity Decision Problem} is NP-complete for chordal graphs (connected).
\end{theorem}

\begin{remark}
	The \textsc{Maximum Transitivity Decision Problem} is known to be NP-complete for the perfect elimination bipartite graphs
	\cite{paul2023transitivity} and doubly chordal graphs \cite{santra2023transitivity}. Therefore, from Lemma\ref{ttr_np_lemma_1}, we can say that \textsc{Maximum Transitivity Decision Problem} is NP-complete for perfect elimination bipartite graphs (disconnected) and doubly chordal graphs (disconnected).
\end{remark}

\section{Tournament transitivity in trees}
In this section, we design a polynomial-time algorithm for finding the tournament transitivity of a given tree $T$. To design the algorithm, we first prove that the tournament transitivity of a tree $T$ can be either $Tr(T)-1$ or $Tr(T)$, where $Tr(T)$ is the transitivity of $T$.

\begin{lemma}\label{tree_tournament_transitivity_lemma_1}
	For a tree $T$, $Tr(T)-1\leq TTr(T)\leq Tr(T)$.
\end{lemma}
\begin{proof}
	According to the tournament transitive partition definition, for any graph G, $TTr(G)\leq Tr(G)$. Hence, for a tree $T$, we have $TTr(T)\leq Tr(T)$. Note that if $Tr(T)=1$, by \cite{hedetniemi2018transitivity} $T$ is a single vertex graph, and in that case, $TTr(T)=1$. Again, when $Tr(T)=2$, $T$ is star graph, according to \cite{hedetniemi2018transitivity}. Let $T$ be a star $S_t$. According to Proposition \ref{Complete_bipartite_graph_ttr}, we know that $TTr(S_t)=1$ if $t=1$, otherwise $TTr(S_t)=2$. Therefore, we assume that $T$ is a tree such that $Tr(T)\geq 3$.

	To prove $TTr(G)\geq Tr(T) -1$ for tree, consider $\pi=\{V_1, V_2, \ldots, V_k\}$ be a transitive partition of $T$ with size $Tr(T)\geq 3$. As $Tr(T)\geq 3$, from \cite{haynes2019transitivity}, we can modify the above transitive partition into another transitive partition, say, $\pi'=\{V_1', V_2', \ldots, V_k'\}$, such that $|V_k'|=|V_{k-1}'|=1$ and $|V_{k-i}'|\leq 2^{i-1}$ for all $2\leq i\leq k-2$. Since $T$ is a tree, $|V_k'|=|V_{k-1}'|=1$, $|V_{k-i}'|=2^{i-1}$ for all $2\leq i\leq k-2$ and $|V_1'|\geq 2^{k-2}$. Let us consider $V_i'$ and $V_j'$ for $i<j\leq k-1$. Let $x\in V_i'$ and $z\in V_{k}'$ such that $xz\in E(T)$. We show that $x$ has no neighbour in $V_j'$. If possible, assume $y_j\in V_j'$ and $xy_j\in E(T)$. As $|V_j'|=2^{k-j-1}=|(V_{j+1}'\cup V_{j+2}'\cup \ldots \cup V_{k}')|$ and $T$ is a tree for $y_j$, there exists $y_{j+s}\in V_{j+s}'$ such that $y_jy_{j+s}\in E(T)$. Now, either $j+s=k$ or $j+s<k$. If $j+s=k$, then $y_{j+s}=z$, and we have a cycle, which is a contradiction. For the other case, when $j+s<k$, similarly we can say that there exists a $y_{j+p}\in V_{j+p}'$ such that $y_{j+s}y_{j+p}\in E(T)$, where $j+p>j+s$. Again, either $j+p=k$ or $j+p<k$. If $j+p=k$, then $y_{j+p}=z$, and we have a cycle, which is a contradiction. If we continue this process, we end up with a cycle, which is a contradiction as $T$ is a tree. So, $x$ has no neighbour in $V_j'$. Therefore, in $\pi'$, $V_j'$ does not dominate $V_i'$ for all $1\leq i<j\leq k-1$. Now consider $\pi''=\{V_1'\cup V_k', V_2', \ldots, V_{k-1}'\}$. Based on the above discussion, we can say that $\pi''$ is a tournament transitive partition of $T$ with size $k-1$. Therefore, $TTr(G)\geq Tr(T)-1$. Hence, for a $T$, $Tr(T)-1\leq TTr(T)\leq Tr(T)$.
\end{proof}


Next, we characterize the trees with tournament transitivity equal to $Tr(T)$. We use a function, \textsc{TransitiveNumber$(x, T)$}, which takes a vertex $x$ and a tree $T$ and returns the transitive number of $x$ in $T$. The transitive number of a vertex $v$ in $T$ is the maximum integer $p$ such that $v\in V_p$ in a transitive partition $\pi=\{V_1, V_2, \ldots, V_k\}$, where the maximum is taken over all transitive partitions of $T$. The transitive number of a vertex $v$ in $T$ is denoted by $t(v, T)$. The function \textsc{TransitiveNumber$(x, T)$} is presented in \cite{hedetniemi1982linear}, where the authors design an algorithm to find the Grundy number of a given tree. As we know that transitivity and Grundy number are the same for a tree \cite{hedetniemi2018transitivity}, the algorithm in \cite{hedetniemi1982linear} also finds the transitivity of a tree. From the description in \cite{hedetniemi1982linear}, it follows that \textsc{TransitiveNumber$(x, T)$} correctly calculates the transitive number of $x$ in $T$.

Let two vertices $u, v\in V(T)$ and the path $uPv$ as $(u, v_{a}, \ldots, v_{b}, v)$. Additionally, we assume that the transitive number of $u$ is $t(u, T)=t(u)$, and the transitive number of $v$ is $t(v, T)=t(v)$. Let $T^{c}$ denote the tree $T$ rooted at $c$ and $T^{[c, c']}$ denote the subtree rooted at $c$, which is obtained by deleting $c'$ from $T^c$. Let us define the set $X$ as the set of vertices required from the path $\{u, v_{a}, \ldots, v_{b}, v\}$ such that $u\in V_{t(u)}$ for some transitive partition $\pi$ of $T$. Clearly, $u\in X$ and the vertices of $X$ are consecutive vertices from the path $(u, v_{a}, \ldots, v_{b}, v)$. We calculate $X$ iteratively by checking whether the transitive number of $u$ is $t(u)-1$ in the tree $T^{[u, v_{\alpha}]}$ or not for all $v_\alpha \in \{v_a, v_{a+1}, \ldots, v_b, v\}$. Also, we define $Y$ as the set of vertices required from $\{u, v_{a}, \ldots, v_{b}, v\}$ to achieve the transitive number of $v$ as $t(v)$ and calculate $Y$ similarly. For a vertex $w\in X\cap Y$, we say that $w$ agrees for $u$ and $v$ if there exist two transitive partitions of $T$ (not necessarily distinct), say $\pi_p$ and $\pi_q$, such that $u\in V_{t(u)}$ and $w\in V_{t}$ in $\pi_p$ and $v\in V_{t(v)}$ and $w\in V_{t}$ in $\pi_q$, respectively. In other words, the index of the sets in which $w$ belongs is the same in both $\pi_p$ and $\pi_q$. In the following lemma, we characterize trees with equal transitivity and tournament transitivity.

\begin{lemma}\label{tree_tournament_transitivity_lemma_2}
	Let $T$ be a tree and $Tr(T)=k(\geq 3)$. Then $TTr(G)=Tr(T)$ if and only if there exist two vertices $y, z$ of $T$ such that the following conditions hold.
	
	\begin{enumerate}
		\item Transitive number of $z$ in $T$ is $k$, that is, $t(z, T)=k$.
		
		\item $y\notin N_T(z)$.
		
		\item Transitive number of $y$ in $T$ is at least $k-1$, that is, $t(y, T)\geq k-1$.
		
		\item Either $X\cap Y=\phi$ or $X\cap Y\neq \phi$ and every vertex of $X\cap Y$ agrees for $y$ and $z$.

		[Where $X$ is the set of vertices required from the path $yPz$ such that $y\in V_{k-1}$ for some transitive partition $\pi$ of $T$ and $Y$ is the set of vertices required from the path $yPz$ such that $z\in V_{k}$ for some transitive partition $\pi$ of $T$. For a vertex $w\in X\cap Y$, we say that $w$ agrees for $y$ and $z$ if there exist two transitive partitions of $T$ (not necessarily distinct), say $\pi_p$ and $\pi_q$, such that $y\in V_{k-1}$ and $w\in V_{t}$ in $\pi_p$ and $z\in V_{k}$ and $w\in V_{t}$ in $\pi_q$, respectively. In other words, the index of the sets in which $w$ belongs is the same in both $\pi_p$ and $\pi_q$.]

	\end{enumerate}
	
\end{lemma}
\begin{proof}
	First, assume that $TTr(T)=Tr(T)=k(\geq 3)$. By the Proposition \ref{ttr_last_two_sets_size}, we have a tournament transitive partition of $T$ of size $k$, say $\pi=\{V_1, V_2, \ldots, V_k\}$, such that $|V_k|=1$ and $|V_{k-1}|=2$. Moreover, $V_{k-1}=\{x, y\}$ and $V_k=\{z\}$ such that $xz\in E(T)$ and $yz\notin E(T)$. Since $\pi$ is also a transitive partition of $T$, by the definition of transitive number, we have $t(z, T)=k$ and $t(y, T)\geq k-1$. Also, by our choice, $y\notin N_T(z)$. Furthermore, $\pi$ is a transitive partition such that $y\in V_{k-1}$ and $z\in V_k$. Now, from the definition of $X$ and $Y$, we have either $X\cap Y=\phi$ or $X\cap Y\neq \phi$, and every vertex of $X\cap Y$ agrees for $y$ and $z$. Therefore, we have two vertices, namely $y$ and $z$, that satisfy the given conditions.
	
	Conversely, assume there exist two vertices $y, z$ of $T$ such that the given conditions hold. Now we divide our proof into the following two cases, based on whether $X\cap Y=\phi$ or $X\cap Y\neq \phi$ and every vertex of $X\cap Y$ agrees for $y$ and $z$. 
	\begin{case}
		$X\cap Y=\phi$
	\end{case}
	Since $X\cap Y=\phi$, we can partition $T=T_1\cup T_2\cup T_3$ such that $T_i\cup T_j=\phi$, for all $i, j$. Moreover, $t(y, T_1)\geq k-1$, and $t(z, T_2)=k$ (see Figure \ref{fig:tree_lemma_case_1}). As $t(y, T_1)\geq k-1$, there exists a transitive partition of $T_1$, say $\{V_1^1, V_2^1, \ldots, V_{k-1}^1\}$, such that $V_{k-1}^1=\{y\}$. Similarly, there exists a transitive partition of $T_2$, say $\{V_1^2, V_2^2, \ldots, V_{k-1}^2, V_k^2\}$, such that $V_{k-1}^2=\{x\}$ and $V_k^2=\{z\}$ for some $x\in V(T_2)$. Now consider the partition $\pi=\{V_1, V_2, \ldots, V_{k-1}, V_k\}$ as follows: $V_1=V(T_3)\cup V_1^1\cup V_1^2$, $V_j=V_j^1\cup V_j^2$ for all $2\leq j\leq k-1$, and $V_k=V_k^2$. By the given conditions and our construction, it is clear that the partition $\pi$ is a tournament transitive partition of $T$. Therefore, in this case $TTr(T)=Tr(T)=k$.

	\begin{figure}[htbp!]
		\centering
		\includegraphics[scale=0.60]{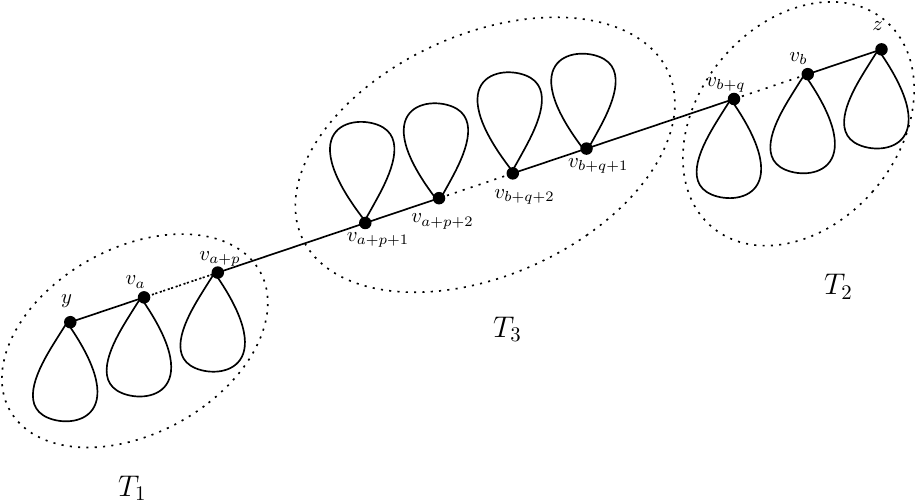}
		\caption{Partition of tree $T$ into $T_1, T_2$ and $T_3$}
		\label{fig:tree_lemma_case_1}
	\end{figure}

	\begin{case}
		$X\cap Y\neq \phi$ and every vertex of $X\cap Y$ agrees for $y$ and $z$. 
	\end{case}
	
	In this case we have two situations when we calculate the transitive number of $z$, that is, $t(z, T)=k$. The neighbour of $z$, which contributes $k-1$ to calculate the transitive number of $z$, is either not in the path $yPz$ or in the path. For the first part, let us assume $x$ is such a neighbour of $z$, which contributes $k-1$ to calculate the transitive number of $z$. Then we can construct a tournament transitive partition $\pi=\{V_1, V_2, \ldots, V_{k-1}, V_k\}$ by setting $V_k=\{x\}$,  $V_{k-1}=\{y, z\}$, and all the other vertices accordingly (see Figure \ref{fig:tree_lemma_case_2_part_1}). Therefore, in this case $TTr(T)=Tr(T)=k$.

	\begin{figure}[htbp!]
		\centering
		\includegraphics[scale=0.65]{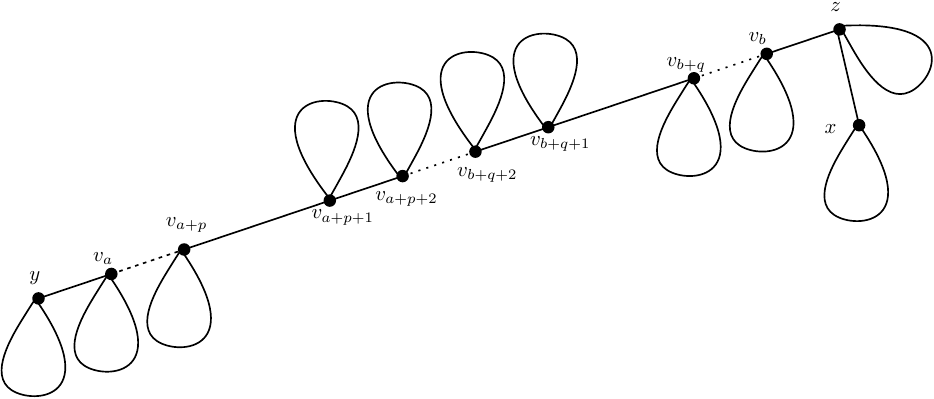}
		\caption{Tournament transitive partition of $T$ when $x$ is not in the path $yPz$}
		\label{fig:tree_lemma_case_2_part_1}
	\end{figure} 
	
	Consider the other part when $x$ is in the path $yPz$. Note that $x$ cannot be a common vertex of $y$ and $z$ as every vertex of $X\cap Y$ agrees for $y$ and $z$. In this situation we can construct a tournament transitive partition $\pi=\{V_1, V_2, \ldots, V_{k-1}, V_k\}$ by setting $V_k=\{z\}$,  $V_{k-1}=\{y, x\}$, and all the other vertices accordingly (see Figure \ref{fig:tree_lemma_case_2_part_2}). 
	
	\begin{figure}[htbp!]
		\centering
		\includegraphics[scale=0.65]{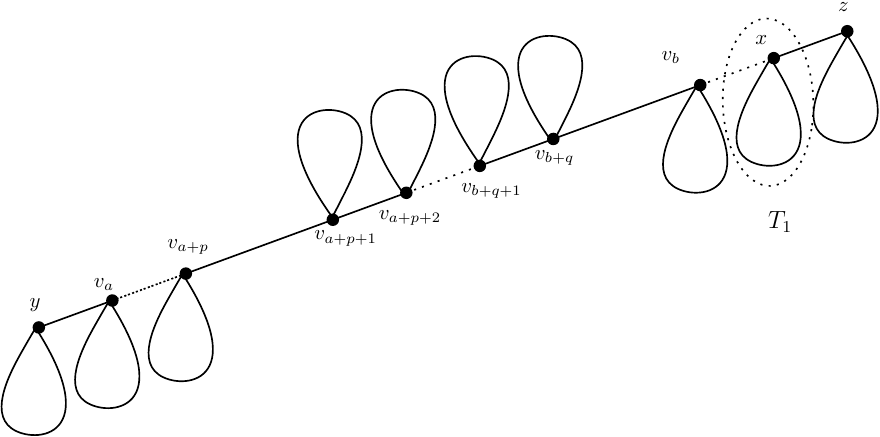}
		\caption{Tournament transitive partition of $T$ when $x$ is in the path $yPz$}
		\label{fig:tree_lemma_case_2_part_2}
	\end{figure} 
	
	From the Figure \ref{fig:tree_lemma_case_2_part_2} it is clear that $z$ has at least one each neighbour in $V_1, V_2, \ldots, V_{k-2}$, which are neither a neighbour of $x$ nor a $y$. Therefore, $V_{k-1}$ does not dominate $V_j$ for all $1\leq j\leq k-2$. Therefore, in this case also $TTr(T)=Tr(T)=k$.

\end{proof}

\subsection{Finding the set $X$ and $Y$}
In this subsection, we calculate the sets $X$ and $Y$, which are defined as above. Now, we investigate those transitive partitions of $T$, where $u\in V_{t(u)}$, and show that for those transitive partitions, the vertices of $X$ can be in some specific sets. To show this result, we need to understand how to find the transitive number of $x$ in a rooted tree $T^x$, and when $x$ achieves that transitive number, then to which sets the children of $x$ belong. The following lemma shows that.
\begin{lemma}\label{tree_lemma_transitivity}
	Let $v_1, v_2, \ldots, v_k$ be the children of $x$ in a rooted tree $T^x$, and for each $1\leq i\leq k$, $l_i$ denotes the rooted transitive number of $v_i$ in $T^x$ with $l_1\leq l_2\leq \ldots\leq l_k$. Let $z$ be the largest integer such that there exists a subsequence of $\{l_i: 1\leq i\leq k\}$, say $(l_{i_1}\leq l_{i_2}\leq \ldots \leq l_{i_{z}})$ such that $l_{i_j}\geq j$, for all $1\leq j\leq z$.
	\begin{enumerate}
		
		\item[(a)] In this case, the transitive number of $x$ in $T^x$ is $1+z$.
		
		\item[(b)] (i) For each $j$ such that $1\leq j\leq z$, there exists a transitive partition $\pi=\{V_1, V_2, \ldots, V_{1+z}\}$ such that $x\in V_{1+z}$ and $v_{i_j}$ belongs to any of the sets from $V_{j}, V_{j+1}, \ldots, V_{\min\{l_{i_j}, z\}}$. 
		
		(ii)  If there exists a neighbour of $x$, say $v$, other than $\{v_{i_1}, v_{i_2}, \ldots, v_{i_z}\}$, such that the rooted transitive number of $v$ in $T^x$ is at least $j$, then there exists a transitive partition $\pi=\{V_1, V_2, \ldots, V_{1+z}\}$ such that $x\in V_{1+z}$ and $v_{i_j}$ belongs to any of the sets from $V_{1}, V_{2}, \ldots, V_{j-1}$. Otherwise, let $r$ be the minimum index such that $l_{i_t}\geq t+1$ for all $r\leq t\leq j-1$. Then there exists a transitive partition $\pi=\{V_1, V_2, \ldots, V_{1+z}\}$ such that $x\in V_{1+z}$ and $v_{i_j}$ belongs to any of the sets from $V_{r}, V_{r+1}, \ldots, V_{j-1}$ if and only if such $r$ exists. 
	\end{enumerate}
	
\end{lemma}

\begin{proof}
	$(a)$ The proof is similar to the result used in \cite{hedetniemi1982linear}, where authors designed an algorithm to calculate the Grundy number of a tree. We rephrase the statement in terms of transitivity.
	
	$(b)(i)$ Let us denote the subtree of $T^x$, rooted at $y$, as $T^x_y$. Since $l_{i_j}\geq j$, for all $1\leq j\leq z$, there exists a transitive partition of $T^x_{v_{i_j}}$, say $\{V_1^j, V_2^j, \ldots, V_{j}^j\}$ such that $v_{i_j}\in V_j^j$. Now let us consider a vertex partition of $T^x$, say $\{V_1, V_2, \ldots, V_{z}, V_{1+z}\}$ as follows: $V_r=\bigcup_{s=r}^zV^s_r$ for all $1\leq r\leq z$ and $V_{1+z}=\{x\}$. Clearly, the above partition is a transitive partition of $T^x$, such that $x\in V_{1+z}$ and $v_{i_j}$ belongs to the set $V_{j}$. Now consider the vertex $v_{i_j}$ and let $s\in \{j+1, j+2, \ldots, \min\{l_{i_j}, z\} \}$. Since $l_{i_j}\geq s$, there exists a transitive partition of $T^x_{v_{i_j}}$, say $\{V_1^j, V_2^j, \ldots, V_{s}^j\}$ such that $v_{i_j}\in V_s^j$. Also, since $l_{i_s}\geq s>j$, there exists a transitive partition of $T^x_{v_{i_s}}$, say $\{V_1^s, V_2^s, \ldots, V_{j}^s\}$ such that $v_{i_s}\in V_j^s$. Similarly, as before, we can show that there exists a transitive partition $\pi=\{V_1, V_2, \ldots, V_{1+z}\}$ such that $x\in V_{1+z}$ and $v_{i_j}$ belongs to the set $V_{s}$. Therefore, we always have a transitive partition of $T^x$, say $\pi=\{V_1, V_2, \ldots, V_{1+z}\}$ such that $x\in V_{1+z}$ and $v_{i_j}$ belongs to any of the sets from $V_{j}, V_{j+1}, \ldots, V_{\min\{l_{i_j}, z\}}$.

	$(b)(ii)$ Let $v$ be a neighbour of $x$ other than $\{v_{i_1}, v_{i_2}, \ldots, v_{i_z}\}$, such that the rooted transitive number of $v$ in $T^x$ is at least $j$. In this case, let $s\in \{1, 2, \ldots, j-1 \}$. Since $l_{i_j}\geq s$, there exists a transitive partition of $T^x_{v_{i_j}}$, say $\{V_1^j, V_2^j, \ldots, V_{s}^j\}$ such that $v_{i_j}\in V_s^j$. Also, since the rooted transitive number of $v$ in $T^x$ is at least $j$, there exists a transitive partition of $T^x_{v}$, say $\{V_1^s, V_2^s, \ldots, V_{j}^s\}$ such that $v\in V_j^s$. By taking the unions of the sets of these partitions as in $(b)(i)$, we can construct a transitive partition of $T^x$, say $\pi=\{V_1, V_2, \ldots, V_{1+z}\}$ such that $x\in V_{1+z}$ and $v_{i_j}$ belongs to the set $V_{s}$.

	Now assume $x$ has no neighbour other than $\{v_{i_1}, v_{i_2}, \ldots, v_{i_z}\}$, having at least $j$ as its rooted transitive number in $T^x$. In this case, let $r$ be the minimum index such that $l_{i_t}\geq t+1$ for all $r\leq t\leq j-1$. For $1\leq s\leq r-1$ and $j+1\leq s\leq z$ as $l_{i_t}\geq t$, we know that there exists a transitive partition of $T^x_{v_{i_s}}$, say $\{V_1^s, V_2^s, \ldots, V_{s}^s\}$ such that $v_{i_s}\in V_s^s$. Again, since $l_{i_t}\geq t+1$ for all $r\leq t\leq j-1$, we have a transitive partition of $T^x_{v_{i_t}}$, say $\{V_1^t, V_2^t, \ldots, V_{t+1}^t\}$ such that $v_{i_t}\in V_{t+1}^t$, for all $r\leq t\leq j-1$. As $l_{i_j}\geq j$, we have a transitive partition of $T^x_{v_{i_j}}$, say $\{V_1^j, V_2^j, \ldots, V_{r}^j\}$ such that $v_{i_j}\in V_r^j$. By taking the unions of the sets of these partitions as in $(b)(i)$, we can construct a transitive partition of $T^x$, say $\pi=\{V_1, V_2, \ldots, V_{1+z}\}$ such that $x\in V_{1+z}$ and $v_{i_j}$ belongs to the set $V_{r}$. Further, using similar arguments as in $(b)(i)$, we can show that there exists a transitive partition $\pi=\{V_1, V_2, \ldots, V_{1+z}\}$ such that $x\in V_{1+z}$ and $v_{i_j}$ belongs to any of the sets from $V_{r+1}, V_{r+2}, \ldots, V_{j-1}$.

	Conversely, let no such $r$ exist. In this case, we show that there does not exist any transitive partition $\pi=\{V_1, V_2, \ldots, V_{1+z}\}$ such that $x\in V_{1+z}$ and $v_{i_j}$ belongs to one of the sets from $V_{1}, V_{2}, \ldots, V_{j-1}$. Let $\pi=\{V_1, V_2, \ldots, V_{1+z}\}$ be a transitive such that $x\in V_{1+z}$ and $v_{i_j}$ belongs to $V_{s}$, for some $1\leq s\leq j-1$. We know that $l_{i_1}\leq l_{i_2}\leq \ldots \leq l_{i_z}$ and $r$ do not exist. This implies that $l_{i_1}\leq l_{i_2}\leq \ldots \leq l_{i_{j-1}}<j$. So, vertices of $\{v_{i_1}, v_{i_2}, \ldots v_{i_{j-1}}\}$ cannot belong to any of the sets $\{V_{j}, V_{j+1}, \ldots, V_z\}$. As we assume that $v_{i_j} \in V_s$ for some $1\leq s\leq j-1$ and all the neighbours of $x$, other than $\{v_{i_1}, v_{i_2}, \ldots, v_{i_z}\}$, have rooted transitive number less than $j$, the vertices $\{v_{i_{j+1}}, v_{i_{j+2}}, \ldots v_{i_{z}}\}$ are the only neighbours of $x$ that can belong to the sets $\{V_{j}, V_{j+1}, \ldots, V_{z}\}$. Clearly, there exists a set from $\{V_{j}, V_{j+1}, \ldots, V_{z}\}$, which does not contain any children of $x$, which contradicts the fact that $\pi$ is a transitive partition. Therefore, there does not exit any transitive partition $\pi=\{V_1, V_2, \ldots, V_{1+z}\}$ such that $x\in V_{1+z}$ and $v_{i_j}$ belongs to any of the sets from $V_{1}, V_{2}, \ldots, V_{j-1}$. Hence, we have the lemma.
\end{proof}

Based on the above lemma, we now show that in a transitive partition of $T$, if $u\in V_{t(u)}$, then the vertices of $X$ can be in some specific sets. A similar result is true for $Y$ as well.

\begin{lemma}\label{properties_of_X_or_Y}
	Let $u$ and $v$ be two vertices of $T$ and they are connected by the path $(u, v_{a}, \ldots, v_{b}, v)$ in $T$. Let $t(u, T)=t(u)$ and $t(v, T)=t(v)$. Further, let $X=\{x_1, x_2, \ldots x_{\alpha}\}$ be the set of consecutive vertices, starting with $x_1=u$, from the path that is required to achieve the transitive number of $u$ is $t(u)$. Similarly, let $Y=\{y_1, y_2, \ldots y_{\beta}\}$ be the set of consecutive vertices, starting with $y_1=v$, from the path that is required to achieve the transitive number of $v$ as $k$. Also, assume that $\pi_p$ and $\pi_q$ are two transitive partitions of $T$ such that $u\in V_{t(u)}$ in $\pi_p$ and $v\in V_{t(v)}$ in $\pi_q$, respectively. Then for all $2\leq j\leq \alpha-1$, $x_j$ belongs to a unique set $V_{p_j}$ in $\pi_p$, where $p_j= t(x_j, T^{[x_j, x_{j-1}]})$ and $x_{\alpha}$ belongs to a set of sets in $\pi_p$ as described in Lemma \ref{tree_lemma_transitivity}(b). Similarly, for all $2\leq j\leq \beta-1$, $y_j$ belongs to a unique set $V_{q_j}$ in $\pi_q$, where $q_j= t(y_j, T^{[y_j, y_{j-1}]})$ and $y_{\beta}$ belongs to a set of sets in $\pi_q$ as described in Lemma \ref{tree_lemma_transitivity}(b).
\end{lemma}

\begin{figure}[h]
	\centering
	\includegraphics[scale=0.70]{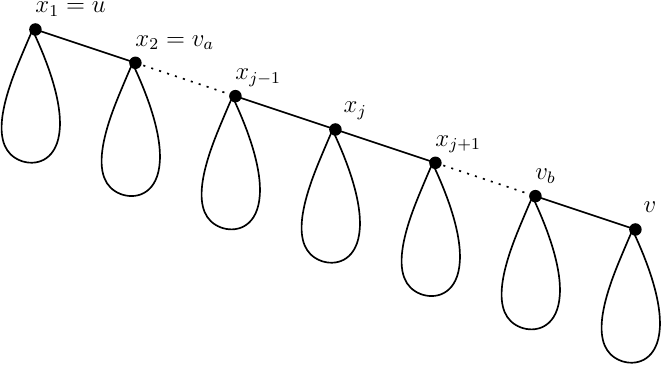}
	\caption{Calculation of X and Y in the tree $T$}
	\label{fig:Calculation_of_X_and_Y}
\end{figure}

\begin{proof}
	Let $x_j\in X$ for $2\leq j\leq \alpha-1$. The transitive number of $x_j$ in $T^{[x_j, x_{j-1}]}$ is given by $t(x_j, T^{[x_j, x_{j-1}]})$. Note that, in $\pi_p$, to obtain the transitive number of $u$ as $t(u)$, the vertex $x_j$ must belong to $V_{p_j}$ in $\pi_p$, where $p_j= t(x_j, T^{[x_j, x_{j-1}]})$ (see Figure \ref{fig:Calculation_of_X_and_Y}). Because if $x_j$ belongs to some $V_r$, where $r< {p_j}$, then $x_{j+1}$ does not belong to $X$ as we have transitive partition of $T^{[x_j, x_{j-1}]}\setminus \{x_{j+1}\}$ such that $x_j\in V_r$. For the vertex, $x_{\alpha}$ can be assigned from a set of sets in $\pi_p$, as described in Lemma \ref{tree_lemma_transitivity}(b). The proof for $Y$ is similar.
\end{proof}

Based on the above two lemmas, we have Algorithm\ref{Algo:tree_X_Y_sets} that describes the process of calculating the sets $X$ and $Y$ for two specific vertices $y$ and $z$ with $t(y, T)\geq k-1$ and $t(z, T)=k$.

\begin{algorithm}[h]
	
	\caption{\textsc{Calculate\_X\_and\_Y}} \label{Algo:tree_X_Y_sets}
	
	\textbf{Input:} The tree $T$ and two vertices  $y$ and $z$ with $t(y, T)\geq k-1$ and $t(z, T)=k$
	
	\textbf{Output:} Set $X$ and $Y$ corresponding $y$ and $z$, respectively
	.
	\begin{algorithmic}[1]

		\State Let the path $P=(v_{a_0}=y, v_{a}, \ldots, v_{b}, z=v_{b_0})$;
		
		
		\State Initially, $X=\{y\}$ and  $Y=\{z\}$;
		
		\If {$t(y)=t(y, T)=k$}
		
		\State $t(y, T^{[y, v_{a}]})\geq k-1$;
		
		\State $X=\{y\}$;
		
		\Else

		\ForAll {$r$ in $\{v_{a}, v_{a+1}, \ldots, v_{b}, z\}$}
		
		\If {$t(y, T^{[y, r]})=k-2$}
		
		\State $X=X\cup \{r\}$;
		
		\Else
		
		\State Stop;
		\EndIf
		
		\EndFor
		
		\EndIf

		\ForAll {$s$ in $\{v_{b}, v_{b-1}, \ldots, v_{a}, y\}$}
		
		\If {$t(z, T^{[z, s]})=k-1$}
		
		\State $Y=Y\cup \{s\}$;
		
		\Else
		
		\State Stop;
		
		\EndIf

		\EndFor
		
		\State  Calculate the set(s) in $\pi_p$ to which the vertices of $X$ belong using to the Lemma \ref{properties_of_X_or_Y};
		
		\State Calculate the set(s) in $\pi_q$ to which the vertices of $Y$ belong using to the Lemma \ref{properties_of_X_or_Y};
	\end{algorithmic}
	
\end{algorithm}

Based on Lemmas \ref{tree_tournament_transitivity_lemma_1} and \ref{tree_tournament_transitivity_lemma_2}, we can design an algorithm for finding tournament transitivity of a given tree. 

\begin{algorithm}[h]
	
	\caption{\textsc{Tournament\_Transitivity(T)}}\label{Algo:tree_tounament_transitivity(T)}
	
	\textbf{Input:} A tree $T$
	
	\textbf{Output:} Tournament transitivity of $T$, that is, $TTr(T)$

	\begin{algorithmic}[1]

		\State Let $V=\{u_1, u_2, \ldots, u_n\}$ be the set of all vertices of $T$;
		
		\ForAll {$u \in V$}
		
		\State \textsc{TransitiveNumber$(u, T)$} ;
		
		\EndFor
		
		\State $S=\{(y, z)| y, z\in V \text{and }y\notin N_T(z) \text { and } t(y, T)\geq k-1, t(z, T)=k \}$
		
		\If {$S\neq \phi$}

		\ForAll {$(y, z)\in S$}

		\State Calculate the sets $X$ and $Y$ according to the Algorithm \ref{Algo:tree_X_Y_sets};
		
		\If {$X\cap Y=\phi$}
		
		\State $TTr(T)=Tr(T)$;
		
		\ElsIf{$X\cap Y\neq \phi$ and every vertex of $X\cap Y$ agrees for $y$ and $z$}
		
		\State $TTr(T)=Tr(T)$;
		
		\Else 
		
		\State $TTr(T)=Tr(T)-1$;
		
		\EndIf
		
		\EndFor
		
		\Else
		
		\State $TTr(T)=Tr(T)-1$;
		
		\EndIf

		\State \Return($TTr(G)$);

	\end{algorithmic}
	
\end{algorithm}

\subsection{Running time of the algorithm}

In this subsection, we analyse the running time of Algorithm \ref{Algo:tree_tounament_transitivity(T)}. In lines $2-4$ of the algorithm, we have computed the transitive number of every vertex. According to the running time analysis of \cite{hedetniemi1982linear}, it will take $O(n)$ to calculate the transitive number of every vertex. Moreover, for a fixed pair of vertices $\{y, z\}$, we can check whether it is in $S$ or not in a constant time.

Now we analyse the running time of Algorithm \ref{Algo:tree_X_Y_sets} (\textsc{Calculate\_X\_and\_Y}) for a pair of vertices $\{y, z\}$ such that $y\notin N_T(z)$ and $t(y, T)\geq k-1, t(z, T)=k$. While computing the set $X$, the checking in line $8$ takes $\displaystyle\sum_{v\in T}(O(deg(v))$ time. This implies that we can compute the set $X$ in $k(\displaystyle\sum_{v\in T}O(deg(v)))$ time, where $k$ is the length of the path $yPz$. In the worst case, $X$ in lines $3-14$ can be computed in $O(n^2)$ time. Similarly, $Y$ can be computed in $O(n^2)$ time. The process of computing $X$ and $Y$ also indicates the set(s) to which the vertices of $X$ and $Y$ belong. Therefore, lines $22$ and $23$ take constant time. So for all vertices from $S$, to calculate the sets $X$ and $Y$, we need at most $n^2*O(n^2)=O(n^4)$ time. As a result, we get that the running time of Algorithm \ref{Algo:tree_tounament_transitivity(T)} is $O(n)+O(n^4)=O(n^4)$. Hence, we have the following main theorem.

\begin{theorem}
	The \textsc{Maximum Tournament Transitivity Problem} can be solved in polynomial time for trees.
\end{theorem}

\section{Tournament transitivity in bipartite chain graphs}

In this section, we characterize some bipartite chain graphs with equal tournament transitivity and transitivity. Let $G=(X\cup Y, E)$  be a bipartite chain graph also let $\sigma_X= (x_1,x_2, \ldots, x_{n_1})$ and $\sigma_Y=(y_1,y_2, \ldots, y_{n_2})$ be the chain ordering. We know from the Proposition \ref{Complete_bipartite_graph_ttr}, for complete bipartite graph $TTr(K_{m, n})=2$, if and only if either $m\neq 1$ or $n\neq 1$. Let us assume only bipartite chain graphs, which are not complete bipartite graphs. Now we find some bipartite chain graphs which are not a complete bipartite graph with equal tournament transitivity and transitivity. Assume $t$ be the maximum integer such that $G$ contains $K_{t, t}$ as an induced subgraph. It is known form \cite{paul2023transitivity}, that $K_{t, t}=G[X_t\cup Y_t]$, where $X_t=\{x_1, x_2, \ldots, x_t\}$ and $Y_t=\{y_1, y_2, \ldots, y_t\}$. Now in $G$, the edges $x_{t+1}y_t, x_ty_{t+1}$ may or may not present. Based on this, a bipartite chain graph, which is not a complete bipartite graph, can be partitioned into three subclasses.

%
%
%

\begin{definition}
Let $G=(X\cup Y, E)$ be a bipartite chain graph which is not a complete bipartite graph. Also, let $\sigma_X= (x_1,x_2, \ldots, x_{n_1})$ and $\sigma_Y=(y_1,y_2, \ldots, y_{n_2})$ be the chain ordering of $G$ and $t$ be the maximum integer such that $G$ contains a $K_{t, t}$ as an induced subgraph. A bipartite chain graph $G$ is called $(i)$ \textup{Type-I BCG} if $x_{t+1}y_t\notin E(G)$ and $x_ty_{t+1}\notin E(G)$,  $(ii)$ \textup{Type-II BCG} if either $x_{t+1}y_t\in E(G)$ or $x_ty_{t+1}\in E(G)$ not both,  $(iii)$ \textup{Type-III BCG} if $x_{t+1}y_t\in E(G)$ and $x_ty_{t+1}\in E(G)$.
\end{definition} 

It is known form \cite{paul2023transitivity}, that for a bipartite chain graph $G$, $Tr(G)=t+1$, where $t$ is the maximum integer such that $G$ contains either $K_{t,t}$ or  $K_{t,t}-\{e\}$ as an induced subgraph. Also, from the definition of tournament transitivity, we know that $TTr(G)\leq Tr(G)$ for a graph $G$. Next, we find the conditions for which the transitivity and tournament transitivity are the same for \textup{Type-I BCG} and for some \textup{Type-II BCG} graphs. Also, we show that for \textup{Type-III BCG}, $TTr(G)<Tr(G)$ always.

\subsection{Tournament transitivity of \textup{Type-I BCG}}

In this subsection, we find the condition under which $TTr(G)=Tr(G)$ for a \textup{Type-I BCG}. For that, we have the following theorem.

\begin{theorem}
Let $G$ be a \textup{Type-I BCG} with $\sigma_X= (x_1,x_2, \ldots, x_{n_1})$ and $\sigma_Y=(y_1,y_2, \ldots, y_{n_2})$ be the chain ordering of $G$. Then $TTr(G)=Tr(G)=t+1$ if and only if there exist vertices $\{z_1, z_2, \ldots, z_t\}$ from $(X\setminus X_t)\cup (Y\setminus Y_t)$ such that $|N(z_1)|=t-1$ and $t-j+1\geq |N(z_j)|\geq t-j$ for all $2\leq j\leq t$.
\end{theorem}

\begin{proof}
First, we assume that there exist vertices $\{z_1, z_2, \ldots, z_t\}$ from the set $\{x_{t+1}, \ldots x_{n_1}, y_{t+1}, \ldots, y_{n_2}\}$ such that $|N(z_1)|=t-1$ and $t-j+1\geq |N(z_j)|\geq t-j$ for all $2\leq j\leq t$. Consider a vertex partition $\pi=\{V_1, V_2, \ldots, V_{t}, V_{t+1}\}$ of $T$ as follows: $\{x_1, y_1, z_t\}\subseteq V_1$; for all $2\leq i\leq t-1$, we set $V_i=\{x_i, y_i, z_{t-i+1}\}$, $V_t=\{x_t, z_1\}$, and $V_t=\{y_t\}$. We put other vertices of $G$ in $V_1$. This partition $\pi$ is illustrated in Figure \ref{fig:BCG_type-I_partition}.

\begin{figure}[htbp!]
	\centering
	\includegraphics[scale=0.80]{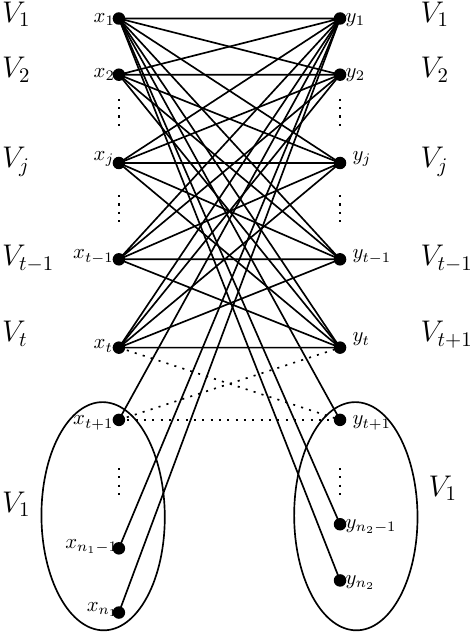}
	\caption{The partition $\pi$ of \textup{Type-I BCG} $G$, where dotted edges are not present in $G$}
	\label{fig:BCG_type-I_partition}
\end{figure} 

We show that $\pi$ is a tournament transitive partition of $G$ with size $t+1$. Since $G$ is a \textup{Type-I BCG}, the set $\{x_{1}, x_{2}, \ldots x_{t}, y_{1}, \ldots, y_{t} \}$ induces a complete bipartite graph $K_{t, t}$ in $G$. Let us consider $V_i$ and $V_j$ for $2\leq i<j\leq t-1$. From $\pi$, we have $V_i=\{x_i, y_i, z_{t-i+1}\}$ and $V_j=\{x_j, y_j, z_{t-j+1}\}$. As $j\geq |N(z_{t-j+1})|\geq j-1$, which implies either $\{x_1, x_2, \ldots, x_{j-1}\}\subseteq N(z_{t-j+1})$ or $\{y_1, y_2, \ldots, y_{j-1}\}\subseteq N(z_{t-j+1})$. So, either $x_i$ is a neighbour of $z_{t-j+1}$ or $y_i$ is a neighbour of $z_{t-j+1}$. Furthermore, because $G$ is a \textup{Type-I BCG}, $x_j$ is adjacent with $y_i$, and $y_j$ is adjacent with $x_i$. Therefore, $V_i$ dominates $V_j$. Moreover, as $i\geq |N(z_{t-i+1})|\geq i-1$, which implies either $N(z_{t-i+1})\subseteq \{x_1, x_2, \ldots, x_{i}\}$ or $N(z_{t-i+1})\subseteq \{y_1, y_2, \ldots, y_{i}\}$. Since \textup{Type-I BCG} and $i<j$, $V_j$ does not dominate $V_i$. 

Now for $V_{t-1}$ and $V_t$, clearly $V_{t-1}$ dominates $V_t$ but $V_t$ does not dominate $V_{t-1}$ as $t\geq |N(z_{t-j+1})|\geq t-1$ and $G$ is a \textup{Type-I BCG}. Similarly, we can show that $V_i$ dominates $V_j$ but $V_j$ does not dominate $V_j$ for others $1\leq i<j\leq k$. Hence, $\pi$ is a tournament transitive partition of $G$ of size $t+1$. Therefore, $TTr(G)\geq t+1$. As we already know from \cite{paul2023transitivity}, for a bipartite chain graph $G$, $Tr(G)=t+1$, where $t$ is the maximum integer such that $G$ contains either $K_{t,t}$ or $K_{t,t}-\{e\}$ as an induced subgraph. Since $G$ is a \textup{Type-I BCG}, $G$ contains $K_{t,t}$ as an induced subgraph for maximum $t$. Hence, $Tr(G)=t+1$. According to the definition of tournament transitivity, we know that $TTr(G)\leq Tr(G)$. Therefore, $TTr(G)=Tr(G)=t+1$.

Conversely, assume $TTr(G)=Tr(G)=t+1$. To show the existence of $\{z_1, z_2, \ldots, z_t\}$ from $\{x_{t+1}, x_{t+2}, x_{n_1}, y_{t+1}, y_{t+2}, \ldots, y_{n_2}\}$ such that $|N(z_1)|=t-1$ and $t-j+1\geq |N(z_j)|\geq t-j$ for all $2\leq j\leq t$, first we prove the following claim.

\begin{claim}\label{BCG_type_I_claim}
	Let $\pi=\{V_1, V_2, \ldots, V_{t}, V_{t+1}\}$ be a tournament transitive partition of $G$ of size $t+1$ such that $|V_{t}|=2$ and $|V_{t+1}|=1$. Then one of the following conditions holds. 
	
	\begin{enumerate}
		\item Each $V_1, \ldots, V_t$ contains exactly one vertex from $\{x_1, x_2, \ldots, x_t\}$ and each $V_1, \ldots, V_{t-1}, V_{t+1}$ contains exactly one vertex from $\{y_1, y_2, \ldots, y_t\}$.
		
		\item Each $V_1, V_2, \ldots, V_t$ contains exactly one vertex from $\{y_1, y_2, \ldots, y_t\}$, and each $V_1, \ldots, V_{t-1}, V_{t+1}$ contains exactly one vertex from $\{x_1, x_2, \ldots, x_t\}$.
	\end{enumerate}
\end{claim}

\begin{proof}
	According to Proposition \ref{ttr_last_two_sets_size}, we always have a $TTr(G)$-partition $\pi=\{V_1, V_2, \ldots, V_{t}, V_{t+1}\}$ of $G$ such that $|V_{t}|=2$ and $|V_{t+1}|=1$. Since $G$ is a \textup{Type-I BCG}, the degree of each vertices from $(X\setminus X_t)\cup (Y\setminus Y_t)$ is at most $t-1$. Now, for any vertex $x\in V_{t+1}$, the degree of $x$ must be at least $t$. Therefore, only vertices from $X_t\cup Y_t$ can be in $V_{t+1}$. Now we prove that if $V_{t+1}$ contains a vertex from $Y_t$, $\pi$ satisfy the condition $(a)$, otherwise $\pi$ satisfy the condition $(b)$.

	Let us  assume $V_{t+1}=\{y_j\}$ for some $1\leq j\leq t$. Since $\pi$ is a tournament transitive partition of $G$, $V_{t}$ must contain a vertex adjacent to $y_j$, and the degree of that vertex is at least $t$. So, for some $1\leq i\leq t$, assume $x_i\in V_{t}$. Moreover, $V_{t+1}$ does not dominate $V_{t}$, which implies that there exists a vertex in $V_{t}$ that is not adjacent to $y_j$, and the degree of that vertex is at least $t-1$. Let $v\in V_{t}$ other than $x_i$, then $v\in (X\setminus X_t)\cup Y$. If $v\in (X\setminus X_t)$, to dominate $v$, each $V_1, V_2, \ldots, V_{t-1}$ contains exactly one vertex from $\{y_1, y_2, \ldots, y_{t-1}\}$, and the vertex $y_j$ must be the vertex $y_t$. As $y_t\in V_{t+1}$ and $N(y_t)=X_t$, each $V_1, V_2, \ldots, V_t$ contains exactly one vertex from $\{x_1, x_2, \ldots, x_t\}$. Therefore, in this case, the condition $(a)$ holds. On the other hand, if $v\in (Y\setminus Y_t)$, we can similarly show the condition $(a)$.

	Finally, assume $v\in Y_t$. As $V_{t-1}$ dominates $V_t$, for $v\in V_t$ there must exist a vertex from $V_{t-1}$, say $u$, such that $v$ is adjacent to $u$ and $|N(u)\setminus \{v, y_j\}|\geq t-2$. If $u\in (X\setminus X_t)$, then $|N(u)\setminus \{v, y_j\}|=|N(u)|-2\leq t-1-2=t-3$, which implies $u$ must be a vertex from $X_t$. Let us assume $u=x_r$ for some $1\leq r\leq t$ and $r\neq i$. As we know, $V_{t}$ does not dominate $V_{t-1}$, so there exists a vertex from $V_{t-1}$, say $v'$, such that $v'$ is not adjacent to the vertices from $\{x_i, v\}$. Therefore, $v'$ must be a vertex from $(X\setminus X_t)\cup (Y\setminus Y_t)$. Consider the case when $v'\in (X\setminus X_t)$. Since $v'\in V_{t-1}$, $|N(v')\setminus \{v, y_j\}|\geq t-2$. But $|N(v')\setminus \{v, y_j\}|=|N(v')|-2\leq t-1-2=t-3$. So, $v'$ cannot be a vertex from $(X\setminus X_t)$. Furthermore, for $v'\in (Y\setminus Y_t)$, $|N(v')\setminus \{x_r, x_i\}|\geq t-2$. But $|N(v')\setminus \{x_r, x_i\}|=|N(v')|-2\leq t-1-2=t-3$. So, $v'$ cannot be a vertex from $(X\setminus X_t)$. We have a contradiction, which implies that $v$ cannot be in $V_t$. Hence, we have the claim.
\end{proof}

From the Claim \ref{BCG_type_I_claim}, we have either each $V_1, V_2, \ldots, V_t$ contains exactly one vertex from $\{x_1, x_2, \ldots, x_t\}$ and each $V_1, V_2, \ldots, V_{t-1}, V_{t+1}$ contains exactly one vertex from $\{y_1, y_2, \ldots, y_t\}$ or each $V_1, V_2, \ldots, V_t$ contains exactly one vertex from $\{y_1, y_2, \ldots, y_t\}$ and each $V_1, V_2, \ldots, V_{t-1}, V_{t+1}$ contains exactly one vertex from $\{x_1, x_2, \ldots, x_t\}$. Since $G$ is a bipartite chain graph, without loss of generality assume $x_r, y_r\in V_r$ for all $1\leq r\leq t-1$ and either $x_t\in V_t$ and $y_t\in V_{t+1}$ or $y_t\in V_t$ and $x_t\in V_{t+1}$. 

Now we are ready to show the existence of $\{z_1, z_2, \ldots, z_t\}$ vertices from $(X\setminus X_t)\cup (Y\setminus Y_t)$ such that $|N(z_1)|=t-1$ and $t-j+1\geq |N(z_j)|\geq t-j$ for all $2\leq j\leq t$. Since $\pi$ is a tournament transitive partition, $V_{t+1}$ does not dominate $V_t$. So, there must exist a vertex from $\{x_{t+1}, x_{t+2}, x_{n_1}, y_{t+1}, y_{t+2}, \ldots, y_{n_2}\}$ that belongs to $V_t$ and is not adjacent to either $y_{t}$ or $x_{t}$, depending on whether $V_{t+1}=\{y_{t}\}$ or $V_{t+1}=\{x_{t}\}$. Let $z_1$ be such a vertex from $(X\setminus X_t)\cup (Y\setminus Y_t)$. As $z_1\in V_{t}$, $|N(z_1)|\geq t-1$ and we know that $|N(z_1)|\leq t-1$, which implies $|N(z_1)|=t-1$. Let us consider $V_j$ and $V_{j+1}$ for some $1\leq j\leq t-1$. Since $V_{j+1}$ does not dominate $V_j$, there exists a vertex form $(X\setminus X_t)\cup (Y\setminus Y_t)$, say $z_{t-j+1}\in V_j$ such that $z_{t-j+1}\notin N(V_{j+1})$. As $z_{t-j+1}\in V_j$ and $\pi$ is a tournament transitive partition of $G$, $|N(z_{t-j+1})|\geq j-1$. Moreover, $z_{t-j+1}\notin N(V_{j+1})$ and $G$ is a bipartite chain graph, implies that $z_{t-j+1}\notin N(x_{p})$ and $z_{t-j+1}\notin N(y_{p})$, for all $j+1\leq p\leq t$. So, $|N(z_{t-j+1})|\leq j$. Therefore, $j\geq |N(z_{t-j+1})|\geq j-1=t-s+1\geq |N(z_{s})|\geq t-s$. Hence, we have $\{z_1, z_2, \ldots, z_t\}$ from $\{x_{t+1}, x_{t+2}, x_{n_1}, y_{t+1}, y_{t+2}, \ldots, y_{n_2}\}$ such that $|N(z_1)|=t-1$ and $t-s+1\geq |N(z_s)|\geq t-s$ for all $2\leq s\leq t$.
\end{proof}

\subsection{Tournament transitivity of \textup{Type-II BCG}}

In this subsection, we find a sufficient condition under which $TTr(G)=Tr(G)$ for \textup{Type-II BCG}. Let $G$ be a \textup{Type-II BCG} with $\sigma_X= (x_1,x_2, \ldots, x_{n_1})$ and $\sigma_Y=(y_1,y_2, \ldots, y_{n_2})$ be the chain ordering of $G$ and $x_ty_{t+1}\in E(G)$. Further, we divide a \textup{Type-II BCG} into two subclasses, based on whether $x_{t+1}y_{t-1}\notin E(G)$ or $x_{t+1}y_{t-1}\in E(G)$. A \textup{Type-II BCG} $G$ is called $(i)$ \textup{Type-II(a) BCG} if $x_{t+1}y_{t-1}\notin E(G)$, $(ii)$ \textup{Type-II(b) BCG} if $x_{t+1}y_{t-1}\in E(G)$. The following theorems find a sufficient condition under which $TTr(G)=Tr(G)$ for a \textup{Type-IIBCG}.

\begin{theorem}
Let $G$ be a \textup{Type-II(a) BCG} with $\sigma_X= (x_1,x_2, \ldots, x_{n_1})$ and $\sigma_Y=(y_1,y_2, \ldots, y_{n_2})$ be the chain ordering of $G$. Then $TTr(G)=Tr(G)=t+1$ if there exist vertices $\{z_1, z_2, \ldots, z_{t-1}\}$ from $(X\setminus X_t)\cup \{y_{t+2}, y_{t+3}, \ldots, y_{n_2}\}$ such that $t-j\geq |N(z_j)|\geq t-j-1$ for all $1\leq j\leq t-1$.
\end{theorem}

\begin{proof}
Let us assume that there exist vertices  $\{z_1, z_2, \ldots, z_{t-1}\}$ from the set $\{x_{t+1}, x_{n_1}, y_{t+2}, \ldots, y_{n_2}\}$ such that $t-j\geq |N(z_j)|\geq t-j-1$ for all $1\leq j\leq t-1$. As $G$ is a \textup{Type-II(a) BCG}, $x_ty_{t+1}\in E(G)$ and $x_{t+1}y_{t-1}\notin E(G)$. Consider a vertex partition $\pi=\{V_1, V_2, \ldots, V_{t}, V_{t+1}\}$ of $G$ as follows: $\{x_1, y_1, z_{t-1}\}\subseteq V_1$; for all $2\leq i\leq t-1$, we set $V_i=\{x_i, y_i, z_{t-i}\}$, $V_t=\{x_t, y_t\}$, and $V_{t+1}=\{y_{t+1}\}$. We put the other vertices of $G$ in $V_1$. This partition $\pi$ is illustrated in Figure \ref{fig:BCG_type-II_a_partition}.

\begin{figure}[htbp!]
	\centering
	\includegraphics[scale=0.80]{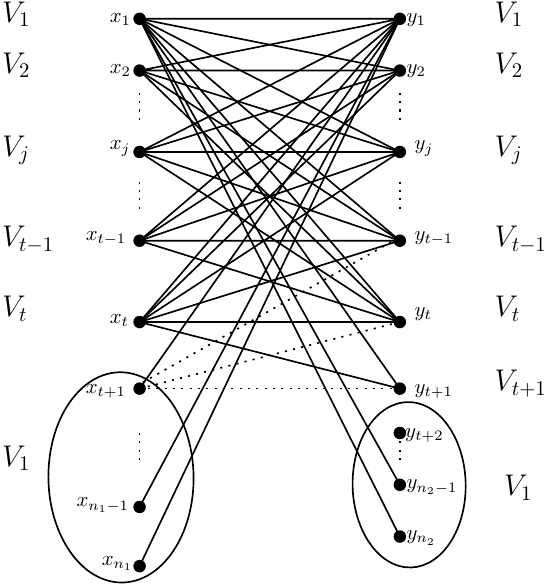}
	\caption{The partition $\pi$ of \textup{Type-II(a) BCG} $G$, where dotted edges are not present in $G$}
	\label{fig:BCG_type-II_a_partition}
\end{figure}

We show that $\pi$ is a tournament transitive partition of $G$ with size $t+1$. Since $G$ is a \textup{Type-II BCG}, the set $\{x_{1}, x_{2}, \ldots x_{t}, y_{1}, \ldots, y_{t} \}$ induces a complete bipartite graph $K_{t, t}$ in $G$. Let us consider $V_i$ and $V_j$ for $2\leq i<j\leq t-1$. From $\pi$, we have $V_i=\{x_i, y_i, z_{t-i}\}$ and $V_j=\{x_j, y_j, z_{t-j}\}$. As $j\geq |N(z_{t-j})|\geq j-1$, which implies either $\{x_1, x_2, \ldots, x_{j-1}\}\subseteq N(z_{t-j})$ or $\{y_1, y_2, \ldots, y_{j-1}\}\subseteq N(z_{t-j})$. Since $i<j$, either $x_i$ is a neighbour of $z_{t-j}$ or $y_i$ is a neighbour of $z_{t-j}$. Furthermore, because $G$ is a \textup{Type-II BCG}, $x_j$ is adjacent with $y_i$, and $y_j$ is adjacent with $x_i$. Therefore, $V_i$ dominates $V_j$ in $\pi$. Moreover, as $i\geq |N(z_{t-i})|\geq i-1$, which implies either $N(z_{t-i})\subseteq \{x_1, x_2, \ldots, x_{i}\}$ or $N(z_{t-i})\subseteq \{y_1, y_2, \ldots, y_{i}\}$. Since \textup{Type-II BCG} and $i<j$, $V_j$ does not dominate $V_i$. Now for $V_{t}$ and $V_{t+1}$, clearly $V_{t}$ dominates $V_{t+1}$ but $V_{t+1}$ does not dominate $V_{t}$. Similarly, we can show that $V_i$ dominates $V_j$ but $V_j$ does not dominate $V_j$ for others $1\leq i<j\leq k$. Hence, $\pi$ is a tournament transitive partition of $G$ of size $t+1$. Therefore, $TTr(G)\geq t+1$. As we already know from \cite{paul2023transitivity}, for a bipartite chain graph $G$, $Tr(G)=t+1$, where $t$ is the maximum integer such that $G$ contains either $K_{t,t}$ or $K_{t,t}-\{e\}$ as an induced subgraph. Since $G$ is a \textup{Type-II BCG}, $G$ contains $K_{t,t}$ as an induced subgraph for maximum $t$. Hence, $Tr(G)=t+1$. According to the definition of tournament transitivity, we know that $TTr(G)\leq Tr(G)$. Therefore, $TTr(G)=Tr(G)=t+1$.
\end{proof}

For \textup{Type-II(b) BCG}, we have the following theorem.

\begin{theorem}
Let $G$ be a \textup{Type-II(b) BCG} with $\sigma_X= (x_1,x_2, \ldots, x_{n_1})$ and $\sigma_Y=(y_1,y_2, \ldots, y_{n_2})$ be the chain ordering of $G$. Then $TTr(G)=Tr(G)=t+1$ if there exist $\{z_1, z_2, \ldots, z_{t-2}\}$ from $\{x_{t+2}, x_{t+3}, \ldots, x_{n_1}, y_{t+2}, y_{t+3}, \ldots, y_{n_2}\}$ such that $t-j-1\geq |N(z_j)|\geq t-j-2$ for all $1\leq j\leq t-2$.
\end{theorem}

\begin{proof}

Let us assume that there exist vertices $\{z_1, z_2, \ldots, z_{t-2}\}$ from the set $\{x_{t+2}, \ldots, x_{n_1}, y_{t+2}, \ldots, y_{n_2}\}$ such that $t-j-1\geq |N(z_j)|\geq t-j-2$ for all $1\leq j\leq t-2$. As $G$ is a \textup{Type-II(b) BCG}, $x_ty_{t+1}\in E(G)$ and $x_{t+1}y_{t-1}\in E(G)$. Consider a vertex partition $\pi=\{V_1, V_2, \ldots, V_{t}, V_{t+1}\}$ of $G$ as follows: $\{x_1, y_1, z_{t-2}\}\subseteq V_1$; for all $2\leq i\leq t-2$, we set $V_i=\{x_i, y_i, z_{t-i-1}\}$, $V_{t-1}=\{x_{t-1}, y_{t-1}\}$, $V_t=\{x_t, x_{t+1}\}$, and $V_{t+1}=\{y_{t}\}$. We put the other vertices of $G$ in $V_1$. This partition $\pi$ is illustrated in Figure \ref{fig:BCG_type-II_b_partition}.

\begin{figure}[htbp!]
	\centering
	\includegraphics[scale=0.80]{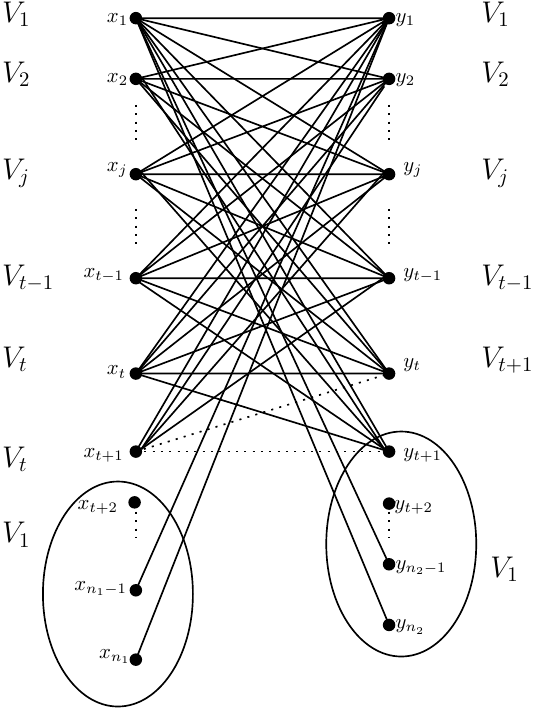}
	\caption{The partition $\pi$ of \textup{Type-II(b) BCG} $G$, where dotted edges are not present in $G$}
	\label{fig:BCG_type-II_b_partition}
\end{figure}

We show that $\pi$ is a tournament transitive partition of $G$ with size $t+1$. Since $G$ is a \textup{Type-II BCG}, the set $\{x_{1}, x_{2}, \ldots x_{t}, y_{1}, \ldots, y_{t} \}$ induces a complete bipartite graph $K_{t, t}$ in $G$. Let us consider $V_i$ and $V_j$ for $2\leq i<j\leq t-2$. From $\pi$, we have $V_i=\{x_i, y_i, z_{t-i-1}\}$ and $V_j=\{x_j, y_j, z_{t-j-1}\}$. As $j\geq |N(z_{t-j-1})|\geq j-1$, which implies either $\{x_1, x_2, \ldots, x_{j-1}\}\subseteq N(z_{t-j-1})$ or $\{y_1, y_2, \ldots, y_{j-1}\}\subseteq N(z_{t-j-1})$. Since $i<j$, either $x_i$ is a neighbour of $z_{t-j-1}$ or $y_i$ is a neighbour of $z_{t-j-1}$. Furthermore, because $G$ is a \textup{Type-II BCG}, $x_j$ is adjacent with $y_i$, and $y_j$ is adjacent with $x_i$. Therefore, $V_i$ dominates $V_j$ in $\pi$. Moreover, as $i\geq |N(z_{t-i-1})|\geq i-1$, which implies either $N(z_{t-i-1})\subseteq \{x_1, x_2, \ldots, x_{i}\}$ or $N(z_{t-i-1})\subseteq \{y_1, y_2, \ldots, y_{i}\}$. Since \textup{Type-II BCG} and $i<j$, $V_j$ does not dominate $V_i$. Now for $t-1\leq p<q\leq t+1$, as $x_{t+1}y_{t-1}\in E(G)$, clearly, $V_p$ dominates $V_q$ and $V_q$ does not dominate $V_p$. Similarly, we can show that $V_i$ dominates $V_j$ but $V_j$ does not dominate $V_j$ for others $1\leq i<j\leq k$. Hence, $\pi$ is a tournament transitive partition of $G$ of size $t+1$. Therefore, $TTr(G)\geq t+1$. As we already know from \cite{paul2023transitivity}, for a bipartite chain graph $G$, $Tr(G)=t+1$, where $t$ is the maximum integer such that $G$ contains either $K_{t,t}$ or $K_{t,t}-\{e\}$ as an induced subgraph. Since $G$ is a \textup{Type-II BCG}, $G$ contains $K_{t,t}$ as an induced subgraph for maximum $t$. Hence, $Tr(G)=t+1$. According to the definition of tournament transitivity, we know that $TTr(G)\leq Tr(G)$. Therefore, $TTr(G)=Tr(G)=t+1$.
\end{proof}

\subsection{Tournament transitivity of \textup{Type-III BCG}}

In this subsection, we show that $TTr(G)<Tr(G)$ for a \textup{Type-III BCG}. For that, we have the following theorem.

\begin{theorem}
Let $G$ be a \textup{Type-III BCG}. Then $TTr(G)<Tr(G)$.
\end{theorem}

\begin{proof}
By the definition of \textup{Type-III BCG}, $G$ contains $K_{t, t}$ for maximum $t$ and $x_{t+1}y_t, x_ty_{t+1} \in E(G)$. From \cite{paul2023transitivity}, as $G$ contains $K_{t+1, t+1}\setminus \{e\}$ for maximum $t$, we have $Tr(G)=t+2$. Now we show that $TTr(G)$ cannot be $t+2$. If possible assume $TTr(G)=t+2$ and let $\pi=\{V_1, V_2, \ldots, V_{t}, V_{t+1}, V_{t+2}\}$ be a tournament transitive partition of $G$. By the Proposition \ref{ttr_last_two_sets_size} we can assume that $|V_{t+1}|=2$ and $|V_{t+2}|=1$.

Since $G$ is a \textup{Type-III BCG}, the degree of each vertices from $(X\setminus X_t)\cup (Y\setminus Y_t)$ is at most $t$. Now, for any vertex $u\in V_{t+2}$, $u$ must have a degree of at least $t+1$. Therefore, only vertices from $X_t\cup Y_t$ can be in $V_{t+2}$ (See Figure \ref{fig:BCG_type-III_partition}).

\begin{figure}[htbp!]
	\centering
	\includegraphics[scale=0.80]{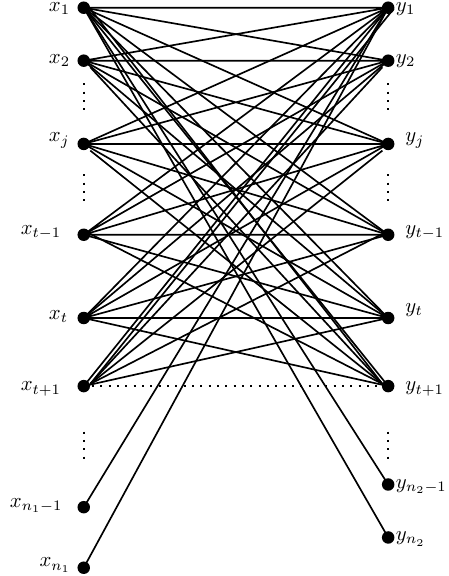}
	\caption{The partition $\pi$ of \textup{Type-III BCG} $G$, where dotted edges are not present in $G$}
	\label{fig:BCG_type-III_partition}
\end{figure}

Let us  assume $V_{t+2}=\{y_j\}$ for some $1\leq j\leq t$. Since $\pi$ is a tournament transitive partition of $G$, $V_{t+1}$ must contain a vertex adjacent to $y_j$, and the degree of that vertex is at least $t+1$. So, for some $1\leq i\leq t$, assume $x_i\in V_{t+1}$. Moreover, $V_{t+2}$ does not dominate $V_{t+1}$, which implies that there exists a vertex in $V_{t+1}$ that is not adjacent to $y_j$, and the degree of that vertex is at least $t$. Let $z\in V_{t+1}$ other than $x_i$, then $z\in (X\setminus X_t)\cup Y$. 

If $z\in (X\setminus X_t)$, to dominate $z$, each $V_1, V_2, \ldots, V_{t}$ contains exactly one vertex from $\{y_1, y_2, \ldots, y_{t}\}$, and the vertex $y_j$ must be the vertex $y_{t+1}$. As $y_{t+1}\in V_{t+2}$ and $N(y_{t+1})=X_t$, each $V_1, V_2, \ldots, V_{t+1}$ contains exactly one vertex from $X_t$. We have a contradiction as $|X_t|=t$, and we require $t+1$ vertices. Assume $z\in (Y\setminus Y_t)$. To dominate $z$, each $V_1, V_2, \ldots, V_{t}$ contains exactly one vertex from $\{x_1, x_2, \ldots, x_{t}\}\setminus \{x_i\}$. Again, we have  contradiction as $|X_t|\setminus \{x_i\}=t-1$ and we required $t$ vertices to dominate $z$.

Finally, assume $z\in Y_t$. As $V_{t}$ dominates $V_{t+1}$, for $z\in V_{t+1}$ there must exist a vertex from $V_{t}$, say $z'$, such that $z$ is adjacent to $z'$ and $|N(z')\setminus \{z, y_j\}|\geq t-1$. If $z'\in (X\setminus X_t)$, then $|N(z')\setminus \{z, y_j\}|=|N(z')|-2\leq t-2=t-2$, which implies $z'$ must be a vertex from $X_t$. Let us assume $z'=x_r$ for some $1\leq r\leq t$ and $r\neq i$. As we know, $V_{t+1}$ does not dominate $V_{t}$, so there exists a vertex from $V_{t}$, say $z''$, such that $z''$ is not adjacent to the vertices from $\{x_i, z\}$. Therefore, $z''$ must be a vertex from $(X\setminus (X_t\cup\{x_{t+1}\}))\cup (Y\setminus (Y_t \cup\{y_{t+1}\}))$. Consider the case when $z''\in (X\setminus (X_t \cup\{x_{t+1}\}))$. Since $z''\in (X\setminus X_t \cup\{x_{t+1}\}))$ and $z''\in V_t$, we must have $|N(z'')\setminus \{z, y_j\}|\geq t-1$. But $|N(z'')\setminus \{z, y_j\}|=|N(z')|-2\leq t-2$. So, $z''$ cannot be a vertex from $(X\setminus (X_t \cup\{x_{t+1}\}))$. Furthermore, for $z''\in (Y\setminus (Y_t \cup\{y_{t+1}\}))$, $|N(z'')\setminus \{x_r, x_i\}|\geq t-1$. But $|N(z'')\setminus \{x_r, x_i\}|=|N(z'')|-2\leq t-2$. So, $z''$ cannot be a vertex from $(Y\setminus (Y_t \cup\{y_{t+1}\}))$. We have a contradiction, which implies that $z$ cannot be in $V_{t+1}$. From the above discussion, we have if $TTr(G)=t+2$, there does not exist any $z\in V_{t+1}$ such that $z$ is not adjacent with $y\in V_{t+1}\cap Y_t$. Similarly, we can show that contradiction when $x\in V_{t+1}\cap X_t$. Hence, $TTr(G)<t+2=Tr(G)$.
\end{proof}

\section{Conclusion and future works}
In this paper, we have studied the notion of tournament transitivity in graphs, which is a variation of transitivity. We have shown that the decision version of this problem is NP-complete for chordal graphs (connected), perfect elimination bipartite graphs (disconnected), and doubly chordal graphs (disconnected). On the positive side, we have proved that this problem can be solved in polynomial time for trees. Furthermore, we have characterized \textup{Type-I BCG} with equal transitivity and tournament transitivity and find some sufficient conditions under which the above two parameters are equal for a \textup{Type-II BCG}. Finally, we have shown that for \textup{Type-III BCG}, these two parameters are never equal.

This paper concludes by addressing some of the several unresolved problems in the study of tournament transitivity of a graph.

\begin{enumerate}
\item What is the necessary condition for a \textup{Type-II BCG} with $TTr(G)=Tr(G)$?

\item We know form \cite{paul2023transitivity}, in linear time we can solved the transitivity problem in bipartite chain graphs. Can we design an algorithm for tournament transitivity in a bipartite chain graph?

\item Characterize connected graphs with $TTr(G)=2$ or $TTr(G)=3$.

\item What is the necessary and sufficient condition for $TTr(G)\geq t$, for an integer $t$?
\end{enumerate}

It would be interesting to investigate the complexity status of this problem in other graph classes. Designing an approximation algorithm for this problem would be another challenging open problem.

\bibliographystyle{plain}
\bibliography{TTR_ref}

\begin{thebibliography}{10}

\bibitem{brandstadt1998dually}
A.~Brandst{\"a}dt, F.~Dragan, V.~Chepoi, and V.~Voloshin.
\newblock Dually chordal graphs.
\newblock {\em SIAM Journal on Discrete Mathematics}, 11(3):437--455, 1998.

\bibitem{chang1994domatic}
G.~J. Chang.
\newblock The domatic number problem.
\newblock {\em Discrete Mathematics}, 125(1-3):115--122, 1994.

\bibitem{CHRISTEN197949}
C.~A. Christen and S.~M. Selkow.
\newblock Some perfect coloring properties of graphs.
\newblock {\em Journal of Combinatorial Theory, Series B}, 27(1):49--59, 1979.

\bibitem{cockayne1977towards}
E.~J. Cockayne and S.~T. Hedetniemi.
\newblock Towards a theory of domination in graphs.
\newblock {\em Networks}, 7(3):247--261, 1977.

\bibitem{effantin2017note}
B.~Effantin.
\newblock A note on grundy colorings of central graphs.
\newblock {\em The Australasian Journal of Combinatorics}, 68(3):346--356,
  2017.

\bibitem{fulkerson1965incidence}
D.~Fulkerson and O.~Gross.
\newblock Incidence matrices and interval graphs.
\newblock {\em Pacific journal of mathematics}, 15(3):835--855, 1965.

\bibitem{furedi2008inequalities}
Z.~F{\"u}redi, A.~Gy{\'a}rf{\'a}s, G.~N. S{\'a}rk{\"o}zy, and S.~Selkow.
\newblock Inequalities for the first-fit chromatic number.
\newblock {\em Journal of Graph Theory}, 59(1):75--88, 2008.

\bibitem{golumbic1978perfect}
M.~C. Golumbic and C.~F. Goss.
\newblock Perfect elimination and chordal bipartite graphs.
\newblock {\em Journal of Graph Theory}, 2(2):155--163, 1978.

\bibitem{haynes2019transitivity}
T.~W. Haynes, J.~T. Hedetniemi, S.~T. Hedetniemi, A.~McRae, and N.~Phillips.
\newblock The transitivity of special graph classes.
\newblock {\em Journal of Combinatorial Mathematics and Combinatorial
  Computing}, 110:181--204, 2019.

\bibitem{haynes2020upper}
T.~W. Haynes, J.~T. Hedetniemi, S.~T. Hedetniemi, A.~McRae, and N.~Phillips.
\newblock The upper domatic number of a graph.
\newblock {\em AKCE International Journal of Graphs and Combinatorics},
  17(1):139--148, 2020.

\bibitem{hedetniemi2018transitivity}
J.~T. Hedetniemi and S.~T. Hedetniemi.
\newblock The transitivity of a graph.
\newblock {\em Journal of Combinatorial Mathematics and Combinatorial
  Computing}, 104:75--91, 2018.

\bibitem{hedetniemi1982linear}
S.~M. Hedetniemi, S.~T. Hedetniemi, and T.~Beyer.
\newblock A linear algorithm for the grundy (coloring) number of a tree.
\newblock {\em Congressus Numerantium}, 36:351--363, 1982.

\bibitem{hedetniemi2004iterated}
S.~M. Hedetniemi, S.~T. Hedetniemi, A.~A. McRae, D.~Parks, and J.~A. Telle.
\newblock Iterated colorings of graphs.
\newblock {\em Discrete Mathematics}, 278(1-3):81--108, 2004.

\bibitem{heggernes2007linear}
P.~Heggernes and D.~Kratsch.
\newblock Linear-time certifying recognition algorithms and forbidden induced
  subgraphs.
\newblock {\em Nordic Journal of Computing}, 14(1-2):87--108, 2007.

\bibitem{paul2023transitivity}
S.~Paul and K.~Santra.
\newblock Transitivity on subclasses of bipartite graphs.
\newblock {\em Journal of Combinatorial Optimization}, 45(1):1--16, 2023.

\bibitem{santra2023transitivity}
S.~Paul and K.~Santra.
\newblock Transitivity on subclasses of chordal graphs.
\newblock In {\em Algorithms and Discrete Applied Mathematics}, pages 391--402,
  Cham, 2023. Springer International Publishing.

\bibitem{samuel2020new}
L.~Samuel and M.~Joseph.
\newblock New results on upper domatic number of graphs.
\newblock {\em Communications in Combinatorics and Optimization},
  5(2):125--137, 2020.

\bibitem{zaker2005grundy}
M.~Zaker.
\newblock Grundy chromatic number of the complement of bipartite graphs.
\newblock {\em The Australasian Journal of Combinatorics}, 31:325--330, 2005.

\bibitem{zaker2006results}
M.~Zaker.
\newblock Results on the grundy chromatic number of graphs.
\newblock {\em Discrete Mathematics}, 306(23):3166--3173, 2006.

\bibitem{zelinka1980domatically}
B.~Zelinka.
\newblock Domatically critical graphs.
\newblock {\em Czechoslovak Mathematical Journal}, 30(3):486--489, 1980.

\bibitem{zelinka1983k}
B.~Zelinka.
\newblock On $ k $-domatic numbers of graphs.
\newblock {\em Czechoslovak Mathematical Journal}, 33(2):309--313, 1983.

\end{thebibliography}

\end{document}